\documentclass[11pt]{amsart}
\usepackage[letterpaper, total={6.5in,9in}, includefoot, centering]{geometry}
\usepackage{amsfonts}
\usepackage{amsmath,amsthm,amssymb,amscd}
\usepackage[foot]{amsaddr}

\usepackage[mathscr]{eucal}
\usepackage{latexsym}
\usepackage[new]{old-arrows}
\usepackage{graphics}
\usepackage{verbatim}
\usepackage{upgreek}
\usepackage{overpic}
\usepackage[all,cmtip]{xy} % comm diagram
\usepackage{enumitem}
\usepackage{tikz}
\usepackage{tikz-cd}
\usepackage[pagebackref]{hyperref}
\renewcommand*\backref[1]{\ifx#1\relax \else (Page #1) \fi}
\usepackage[font=scriptsize]{caption}
\usepackage{upgreek}
\usepackage{algorithm}
\usepackage{algpseudocode}
\makeatletter
\newsavebox{\@brx}
\newcommand{\llangle}[1][]{\savebox{\@brx}{\(\m@th{#1\langle}\)}%
  \mathopen{\copy\@brx\kern-0.5\wd\@brx\usebox{\@brx}}}
\newcommand{\rrangle}[1][]{\savebox{\@brx}{\(\m@th{#1\rangle}\)}%
  \mathclose{\copy\@brx\kern-0.5\wd\@brx\usebox{\@brx}}}
\makeatother
\graphicspath{ {./images/} }

\usepackage[authoryear,sort&compress]{natbib}
\usepackage{adjustbox}
\bibpunct{[}{]}{;}{n}{,}{,}
\usepackage{cancel}
\usepackage{float}
\usepackage{cleveref}

\makeatletter \@addtoreset{equation}{section}

\DeclareMathOperator{\ext}{ext}

\makeatother
\newtheorem{theorem}{Theorem}[section]

\newtheorem{definition}{Definition}[section]
\newtheorem{remark}{Remark}[section]
\newtheorem{example}{Example}[section]
\newtheorem{prop}{Proposition}[section]
\newtheorem{assumption}{Assumption}[section]

\graphicspath{ {./images/} }

\usepackage{xcolor}

\begin{document}
\title{Variational Principles for Hamiltonian Systems}
\author{Brian K. Tran$^\dagger$ and Melvin Leok$^\ddagger$}
\address{$^\dagger$Los Alamos National Laboratory, Theoretical Division, Los Alamos, NM 87545, USA. \newline \indent $^\ddagger$Department of Mathematics, University of California, San Diego, 9500 Gilman Drive, La Jolla, \newline \indent\ \ CA 92093-0112, USA.}
\email{btran@lanl.gov, mleok@ucsd.edu}
\allowdisplaybreaks

\begin{abstract}
Motivated by recent developments in Hamiltonian variational principles, Hamiltonian variational integrators, and their applications such as to optimization and control, we present a new Type II variational approach for Hamiltonian systems, based on a virtual work principle that enforces the Type II boundary conditions through a combination of essential and natural boundary conditions; particularly, this approach allows us to define this variational principle intrinsically on manifolds. We first develop this variational principle on vector spaces and subsequently extend it to parallelizable manifolds, general manifolds, as well as to the infinite-dimensional setting. Furthermore, we provide a review of variational principles for Hamiltonian systems in various settings as well as their applications.
\end{abstract}

\maketitle

\tableofcontents

\section{Introduction}\label{sec:intro}
A major motivation for developing variational principles is the construction of variational integrators. The construction of variational integrators for non-degenerate Hamiltonian systems was traditionally approached from a Lagrangian perspective and applied to a non-degenerate Hamiltonian system via the Legendre transform. Lagrangian variational integrators span a broad class of methods and applications, see, for example, \cite{MaWe2001, Ob2016, deLe2007, FaLeLeMcSc2006, LeLeMc2005a, LeLeMc2007b, LeMaOr2008, LeMaOrWe2004, HaLu2020, XuTrLe2022, HaLe2015, ObSa2015, ZaNeRoSc2019, Stern2015, ZhQiBuBh2014}.

However, many Hamiltonian systems of interest are degenerate and, as such, it is much more natural to construct methods directly from variational principles on the Hamiltonian side and from generating functions based on such variational principles. Furthermore, as observed in \cite{ScLe2017}, even for a non-degenerate Hamiltonian, the method obtained from a Lagrangian variational integrator with a given choice of underlying approximation scheme followed by the discrete Legendre transform may not be equivalent to a Hamiltonian variational integrator with the same underlying approximation scheme due to numerical conditioning issues.

The Type II variational approach to Hamiltonian mechanics was introduced in \cite{LeZh2009}, where they introduce the Type II variational principle on phase space and the associated exact discrete Hamiltonian, which is a generating function of the exact solution with Type II boundary conditions fixing $q(0)$ and $p(T)$ and satisfies the Type II Hamilton--Jacobi equation. In \cite{LeZh2009}, this Type II framework is used to develop a systematic procedure for constructing Hamiltonian variational integrators by approximating the exact discrete Hamiltonian. Crucially, this construction of Hamiltonian variational integrators requires no Lagrangian equivalent and makes no assumptions of non-degeneracy. The properties of such integrators were subsequently analyzed in \cite{ScLe2017}.

The Type II, and similarly Type III, approach to Hamiltonian mechanics has since led to several developments for the study of Hamiltonian systems, particularly those with degenerate Hamiltonians, and in the numerical integration of Hamiltonian systems. Type II Hamiltonian variational principles are developed in \cite{LeZh2009} and their properties are explored in \cite{ScLe2017}. Hamiltonian Taylor variational integrators are introduced in \cite{ScShLe2018}. Hamiltonian variational integrators for stochastic systems are studied in \cite{KrTy2020, HoTy2018}. A Type II variational principle for multisymplectic Hamiltonian PDEs is introduced in \cite{VaLiLe2011} which was used to construct multisymplectic Hamiltonian variational integrators in \cite{TrLe2022}. This Type II perspective has been applied to the study of accelerated optimization, yielding symplectic accelerated optimization methods \cite{SAO1,SAO2,SAO4,SAO5,SAO6}. The Type II perspective is also the natural setting to study adjoint systems arising in adjoint sensitivity analysis and optimal control \cite{TrLe2024adj, TrSoLe2024}.

This paper is organized as follows. In the next section, \Cref{sec:mechanics}, we give a brief review of Lagrangian and Hamiltonian mechanics. Subsequently, the paper is divided into two halves.

In the first half of the paper, \Cref{sec:hamiltonian-type-ii}, we consider Hamiltonian mechanics from the Type II and III perspective and present several new variational principles. In \Cref{sec:boundary-conditions}, we provide a discussion of four types of boundary conditions for Hamiltonian systems (Types I, II, III and IV) as well as the initial value problem. By studying a model degenerate Hamiltonian, we will see that Type I and IV boundary conditions are incompatible with degeneracy whereas Types II and III boundary conditions are compatible. This provides motivation to consider variational principles based on Type II and III boundary conditions, which we consider in \Cref{sec:type-ii-iii-vp}; note that since the Type III case can formally be obtained from the Type II case under time reversal, we focus on the Type II case. In \Cref{sec:type-ii-vp-vector-space}, we first recall the Type II variational principle on vector spaces of \citet{LeZh2009}. As we will see, this variational principle does not make intrinsic sense on a manifold. To remedy these issues, we introduce an alternative Type II variational principle on vector spaces in \Cref{sec:type-ii-vp-d'Alembert} based on a virtual work principle. Subsequently, we show how this new variational principle can be extended to parallelizable manifolds in \Cref{sec:type-ii-vp-parallelizable} and to arbitrary manifolds in \Cref{sec:free-boundary-type-ii}. Finally, we develop an infinite-dimensional version of this new variational principle and of the Type II variational principle of \citet{LeZh2009} in \Cref{sec:inf-dim-type-ii-vp}. 

In the second half of the paper, \Cref{sec:generalizations}, we provide a review of the literature, giving an overview of applications of Hamiltonian variational principles and generalizations to other settings, as well as discussing some connections to the variational principles developed in the first half of the paper. Particularly, we discuss variational integrators in \Cref{sec:variational-integrators}, optimization and control in \Cref{sec:optimization-control}, Hamiltonian PDEs in \Cref{sec:Hamiltonian-PDEs}, constrained Hamiltonian systems in \Cref{sec:constrained-hamiltonian-systems}, and stochastic Hamiltonian systems in \Cref{sec:stochastic}. The review is not intended to be exhaustive but provides several motivating applications for considering Hamiltonian systems from a variational perspective. We refer the reader to the provided references in these sections for more detailed discussion of these topics.

\subsection{Lagrangian and Hamiltonian Mechanics}\label{sec:mechanics}
We first review some basic ideas regarding variational principles in Lagrangian and Hamiltonian mechanics. For a more thorough discussion, see the classic texts \cite{RaMa1987, MaRa1999, Arnold} and for a more recent discussion, see the text \cite{LeLeMc2017} which discusses variational formulations of Lagrangian and Hamiltonian mechanics, with a particular emphasis on global descriptions of the dynamics. Here, we will state the traditional Hamilton's variational principle and Hamilton's phase space variational principle precisely, as we will make analogous statements for the Type II variational principles presented in \Cref{sec:hamiltonian-type-ii}.  Note that we explicitly consider the case of time-dependent Lagrangians and Hamiltonians in this discussion.

Consider a mechanical system with configuration space $Q$, which here we take to be a finite-dimensional vector space. In the Lagrangian description of mechanics on a time interval $[0,T]$, one considers a Lagrangian $L$, which is a mapping
$$ L: [0,T] \times TQ \rightarrow \mathbb{R}. $$
Associated with such a Lagrangian, define the action
\begin{equation}\label{eq:lagrangian-action}
     S_L[q] := \int_0^T L(t,q(t),\dot{q}(t))\, dt.
\end{equation}
The variational approach to derive the equations of motion is to consider curves $q(t)\in Q$ with fixed position endpoints $q(0)=q_0$, $q(T)=q_1$, and demand that the action \eqref{eq:lagrangian-action} is stationary under variations that fix these endpoints. 

More precisely, we define such variations as follows. For a curve $q:[0,T] \rightarrow Q$ satisfying the fixed position endpoint conditions $q(0)=q_0$, $q(T)=q_1$, consider a one-parameter family of curves $q_\epsilon: [0,T] \rightarrow Q$ which preserves these endpoints, i.e.,
\begin{subequations}\label{eq:type-i-one-parameter}
\begin{alignat}{3}
    q_0(t) &= q(t),&\,&&& \text{ for all } t \in [0,T], \\
    q_\epsilon(0) &= q_0,&\, q_\epsilon(T) &= q_1,&\,& \text{ for all } \epsilon.
\end{alignat}
\end{subequations}
The variation associated with this one-parameter family of curves is given by differentiating with respect to $\epsilon$ at $\epsilon = 0$, i.e., for all $t \in [0,T]$,
$$ \delta q(t) := \frac{d}{d\epsilon}\Big|_{\epsilon = 0} q_\epsilon (t).$$
By the properties assumed of the one-parameter family \eqref{eq:type-i-one-parameter}, the variation satisfies $\delta q(0) = 0 = \delta q(T)$ and $\delta q(t)$ is furthermore tangent to $q(t)$. We refer to the set of all such $\delta q$ as Type I variations (Type I refers to the fixed endpoint conditions $q(0)=q_0$, $q(T)=q_1$, as we will discuss further in \Cref{sec:boundary-conditions}). 

For the variational principle, we will consider $C^2$ curves on $Q$ with fixed position endpoints,
$$ C^2_{q_0,q_1}([0,T],Q) := \{ q \in C^2([0,T], Q): q(0)=q_0, q(T)=q_1\}. $$
Type I variations are then elements of the tangent space to this space of curves, 
\begin{align*}
    TC^2_{q_0,q_1}([0,T],Q) =& C^2_{q_0,q_1}([0,T],Q) \\
    & \quad \times \{ \delta q \in C^2([0,T], TQ): \delta q(t) \in T_{q(t)}Q, \delta q(0)=0=\delta q(T)\}.
\end{align*}
A Type I variation of a curve $q$ is given by a tangent vector $\delta q \in T_q C^2_{q_0,q_1}([0,T],Q).$ The flow of such a tangent vector $\delta q$ generates a one-parameter family of curves $q_\epsilon$ satisfying \eqref{eq:type-i-one-parameter}. Thus, associated with this variation, we define the variation of the action $S_L:C^2_{q_0,q_1}([0,T],Q) \rightarrow \mathbb{R}$ whose expression is given by \eqref{eq:lagrangian-action}, to be
$$ \delta S_L[q] := \frac{d}{d\epsilon}\Big|_{\epsilon = 0} S_L[q_\epsilon]. $$
Hamilton's variational principle is then to find a curve $q \in C^2_{q_0,q_1}([0,T],Q)$ such that the action is stationary $\delta S_L[q]=0$ for all Type I variations $\delta q \in T_qC^2_{q_0,q_1}([0,T],Q)$.

\begin{prop}
    Hamilton's variational principle holds if and only if the Euler--Lagrange equation with Type I boundary conditions holds:
    \begin{subequations}\label{eq:euler-lagrange}
    \begin{align}
        \frac{d}{dt} \left( \frac{\partial L}{\partial \dot{q}} \right) - \frac{\partial L}{\partial q} &= 0, \\
        q(0) &= q_0, \\
        q(T) &= q_1.
    \end{align}
    \end{subequations}
\end{prop}
The proof is standard \cite{RaMa1987, MaRa1999, Arnold, LeLeMc2017}, so we omit it.

In the Hamiltonian description of mechanics, we instead consider a Hamiltonian $H$, which is a mapping
$$ H: [0,T] \times T^*Q \rightarrow \mathbb{R}.$$
In the Hamiltonian description, we consider $C^1$ curves $(q,p)$ on $T^*Q$, again with fixed position endpoints $q(0) = q_0$ and $q(T) = q_1$. The relevant space of curves is then
$$ C^1_{q_0,q_1}([0,T],T^*Q) := \{ (q,p) \in C^1([0,T],T^*Q): q(0)=q_0, q(T)=q_1\}. $$
As before, Type I variations are tangent vectors to this space of curves,
\begin{align*}
    TC^1_{q_0,q_1}&([0,T],T^*Q) \\
    &= C^1_{q_0,q_1}([0,T],T^*Q) \times \{ (\delta q,\delta p) \in C^1( [0,T] ,TT^*Q): \\
    & \qquad \qquad \qquad (\delta q(t),\delta p(t)) \in T_{(q(t),p(t))}T^*Q, \delta q(0)= 0 = \delta q(T)\},
\end{align*}
We define the action $S: C^1_{q_0,q_1}([0,T],T^*Q) \rightarrow \mathbb{R}$ by
\begin{equation}\label{eq:classical-hamiltonian-action}
    S[q,p] := \int_0^T [\langle p(t), \dot{q}(t)\rangle - H(t,q(t),p(t))]\,dt,
\end{equation}
\sloppy where $\langle\cdot,\cdot\rangle$ denotes the duality pairing. A variation $(\delta q,\delta p) \in T_{(q,p)}C^1_{q_0,q_1}([0,T],T^*Q)$ generates a one-parameter family of curves $(q_\epsilon,p_\epsilon)$ preserving the fixed position endpoints. For such a variation, we define the variation of the action to be
$$ \delta S[q,p] := \frac{d}{d\epsilon}\Big|_{\epsilon = 0} S[q_\epsilon, p_\epsilon]. $$
Hamilton's phase space variational principle, analogous to before, is the stationarity of the action $\delta S[q,p]=0$ for all Type I variations $(\delta q, \delta p) \in T_{(q,p)}C^1_{q_0,q_1}([0,T],T^*Q)$. 
\begin{prop}
    Hamilton's phase space variational principle holds if and only if Hamilton's equations with Type I boundary conditions hold:
    \begin{subequations}\label{eq:hamilton-equations-type-i}
        \begin{alignat}{2}
            \dot{q}(t) &= D_p H(t,q,p),&\, &\text{ for all } t\in(0,T),\\ 
            \dot{p}(t) &= -D_q H(t,q,p),&\, &\text{ for all } t\in(0,T), \\
            q(0) &= q_0, \\
            q(T) &= q_1.
        \end{alignat}
    \end{subequations}
\end{prop}
Again, the proof is standard so we omit it.

As we will see, for certain Hamiltonians, particularly degenerate ones, the set of curves satisfying equations \eqref{eq:hamilton-equations-type-i} may be the empty set; in other words, there are cases where there are no critical points of the action $S$ with respect to Type I variations. If one is only interested in the equations of motion, this is not an issue, since one can use the variational principle to formally obtain the equations of motion and subsequently consider other boundary conditions which do make sense for the given problem. However, as we discussed in \Cref{sec:intro}, the construction of Hamiltonian variational integrators utilizes variational principles which directly incorporate the boundary conditions. As such, incorporating the boundary conditions directly into the variational principle is useful from this perspective. Deriving such variational principles, particularly for Type II and Type III boundary conditions, will be the main focus of the next section, \Cref{sec:hamiltonian-type-ii}.

\section{Hamiltonian Mechanics from the Type II/III Perspective}\label{sec:hamiltonian-type-ii}
Throughout, we will denote by $Q$ a finite-dimensional vector space and by $M$ a finite-dimensional smooth manifold. Furthermore, for the tangent and cotangent bundles of $Q$, we will use the natural identifications $TQ \cong Q \times Q$ and $T^*Q \cong Q \times Q^*$ and denote the duality pairing by $\langle\cdot,\cdot\rangle: Q^*\times Q \rightarrow \mathbb{R}$. Similarly, we denote the duality pairing between fibers of the tangent and cotangent bundles of $M$ by $\langle\cdot,\cdot\rangle: T^*_qM \times T_qM \rightarrow \mathbb{R}$. We denote by $(q,p)$ canonical coordinates on $T^*M$ or $T^*Q \cong Q \times Q^*$. For a tangent vector $(q,v) \in TQ$ and a cotangent vector $(q,p) \in T^*Q$, we will often use the simpler notation $v \in TQ$ and $p \in T^*Q$, when the base point is clear (and similarly for $TM$, $T^*M$). All maps are assumed to be smooth unless stated otherwise.

\subsection{Boundary Conditions}\label{sec:boundary-conditions}
Throughout this work, we will generally focus on Type II/III variational principles, although we additionally provide an overview of more general variational principles in \Cref{sec:generalizations}. As motivation for such variational principles, we consider the following five types of boundary conditions for a Hamiltonian system on $T^*Q$ over an interval $[0,T]$, as given in \Cref{table:boundary-conditions}.
\begin{table}[H]
\begin{center}
\begin{tabular}{c|c}
     Type 0 & $q(0) = q_0, p(0) = p_0$ \\ \hline
     Type I & $q(0) = q_0, q(T) = q_1$ \\ \hline
     Type II & $q(0) = q_0, p(T) = p_1$ \\ \hline
     Type III & $p(0) = p_0, q(T) = q_1$ \\ \hline
     Type IV & $p(0) = p_0, p(T) = p_1$ 
\end{tabular}
\caption{Five types of boundary conditions for a Hamiltonian system.}
\label{table:boundary-conditions}
\end{center}
\end{table}
In principle, one can consider more general types of boundary conditions (e.g., mixed Robin boundary conditions); however, these five types will cover all of the cases of interest for our purposes. The terminology for such boundary conditions was coined by \citet{goldstein:mechanics} and we adhere to the terminology here.

For the various types of boundary conditions in \Cref{table:boundary-conditions}, we will study their completeness, in the sense of whether they are sufficient to specify a solution curve $(q(t),p(t))$ on the interval $[0,T]$ with $T>0$. In particular, we will study the completeness of these boundary conditions for degenerate Hamiltonian systems, which we define below. The general theory regarding the existence and uniqueness of the solution of two-point boundary value problems is complicated and is not the focus of our discussion (see \cite{BaSh1969, LaOp1967} for early results on existence and uniqueness of two-point boundary value problems and \cite{ElHe2016} for a recent compilation of results); as such, we will consider a prototypical model for a degenerate Hamiltonian system in our analysis and furthermore, we will assume that, for a given regular Hamiltonian system, all five types of boundary conditions are complete on some interval $[0,T]$.

\begin{definition}[Regularity and Degeneracy]
    We say that a Hamiltonian $H: [0,T]\times T^*Q \rightarrow \mathbb{R}$ is regular, or non-degenerate, if for all $t \in [0,T]$ and $(q,p) \in T^*Q$, its Hessian with respect to the momentum coordinates,
    \begin{equation}\label{eq:Hessian-p}
        D_p^2 H(t,q,p),
    \end{equation}
    is invertible. Otherwise, we say that $H$ is degenerate.

    We say that $H$ is maximally degenerate if its Hessian in the momentum coordinates, \Cref{eq:Hessian-p}, is identically zero.

    Furthermore, we say that $H$ is hyperregular if the associated Legendre transform,
    $$ \mathbb{F}H: T^*Q \rightarrow TQ, $$
    defined by the fiber derivative of $H$, $\mathbb{F}H: (q,p) \in T^*_qQ \mapsto D_p H(q,p) \in T_qQ$, is a diffeomorphism at each $t \in [0,T]$. See \cite{RaMa1987, MaRa1999}. 
\end{definition}

\begin{assumption}\label{assumption:complete}
    For a regular Hamiltonian $H$, we will assume that all five types of boundary are complete, i.e., for the system given by Hamilton's equation and any one of the five types of boundary conditions, there exists a unique solution $(q(t),p(t))$ on $[0,T]$. Note that for $T$ sufficiently small, this is generally true for all five types of boundary conditions, by using an implicit function theorem argument with the associated generating functions.
\end{assumption}

Note that a maximally degenerate Hamiltonian must be of the form
$$ H(t,q,p) = \langle p,f(t,q)\rangle + g(t,q), $$
for some functions $f: [0,T] \times Q \rightarrow Q$ and $g: [0,T]\times Q\rightarrow \mathbb{R}$, which follows from integrating the equation $D_p^2 H \equiv 0$ twice with respect to $p$.

We now define our model problem for a degenerate Hamiltonian system, essentially to be the sum of a regular term and a maximally degenerate term. More precisely, let $Q_r, Q_d$ be two vector spaces with $\dim(Q_r) > 0, \dim(Q_d) > 0$ and let $Q = Q_r \times Q_d$; then, $T^*Q \cong T^*Q_r \times T^*Q_d$, with coordinates $(q,p)=(q_r,q_d,p_r,p_d) \in T^*Q$ where $(q_r,p_r) \in T^*Q_r, (q_d,p_d) \in T^*Q_d$. Our model Hamiltonian $H: [0,T] \times T^*Q \rightarrow \mathbb{R}$ is defined by
\begin{align}\label{eq:model-hamiltonian}
    H(t,q,p) = H_r(t,q_r,p_r) + H_d(t,q_d,p_d),
\end{align}
where $H_r$ is a regular Hamiltonian on $T^*Q_r$ and $H_d$ is a maximally degenerate Hamiltonian on $T^*Q_d$, which, as noted above, has the form
$$ H_d(t,q_d,p_d) = \langle p_d, f(t,q_d)\rangle + g(t,q_d). $$
Note that $H$ is degenerate but not maximally degenerate, as seen by computing the Hessian in the momentum coordinates, expressed in block form as
\begin{equation*}
    D_p^2 H(t,q,p) = \begin{pmatrix} D_{p_r}^2 H_r(t,q_r,p_r) & 0 \\ 0 & 0 \end{pmatrix},
\end{equation*}
which has rank $\dim(Q_r)>0$ since $H_r$ is regular. Hamilton's equations for the model Hamiltonian \eqref{eq:model-hamiltonian} are
\begin{subequations}\label{eq:model-ham-eom}
\begin{align}
    \dot{q}_r &= D_{p_r} H_r(t,q_r,p_r), \label{eq:model-ham-eom-a}\\
    \dot{p}_r &= -D_{q_r} H_r(t,q_r,p_r), \label{eq:model-ham-eom-b}\\
    \dot{q}_d &= D_{p_d} H_d(t,q_d,p_d) = f(t,q_d),\label{eq:model-ham-eom-c} \\
    \dot{p}_d &= -D_{q_d} H_d(t,q_d,p_d) = - [D_{q_d}f(t,q_d)]^*p_d - D_{q_d} g(t,q_d), \label{eq:model-ham-eom-d}
\end{align}
\end{subequations}
where $A^*$ denotes the adjoint of $A$ with respect to the duality pairing $\langle\cdot,\cdot\rangle$, i.e., $\langle p, A v \rangle = \langle A^*p,v\rangle$.

We now study the completeness of the five types of boundary conditions for the model Hamiltonian \eqref{eq:model-hamiltonian}.

\textbf{Type 0.} The case of Type 0 boundary conditions is simply the initial value problem $\dot{q} = \partial H/\partial p, \dot{p} = -\partial H/\partial q$ with $q(0)=q_0, p(0)=p_0$. This is the simplest case; indeed, the existence and uniqueness theory for two-point boundary value problems relies on the existence and uniqueness theory for the associated initial value problem \cite{BaSh1969, LaOp1967, ElHe2016}. In our case, for any sufficiently smooth Hamiltonian (not restricted to the above model Hamiltonian), there exists a unique solution curve $(q(t),p(t))$ on $[0,T]$ corresponding to the initial value problem, by the standard existence and uniqueness theory for initial value problems. Thus, Type 0 boundary conditions are complete, in the sense that they are sufficient to specify a unique solution curve to Hamilton's equations on $[0,T]$.

\begin{remark}
    For the remaining cases, we will assume that for any initial value problem encountered below, $\dot{x} = h(t,x), x(0) = x_0$, has a unique solution on $[0,T]$, i.e., that $h$ is sufficiently smooth.
\end{remark}

\textbf{Type I.} Type I boundary conditions are given by specifying the position at the initial and final time, $q(0) = q_0, q(T) = q_1$. We will show that these boundary conditions are incomplete, i.e., insufficient to specify a solution curve for the model Hamiltonian. To see this, note that specifying $q(0) = q_0$ and $q(T) = q_1$ is equivalent to specifying $q_r(0), q_d(0), q_r(T), q_d(T)$. From \Cref{assumption:complete}, the boundary conditions $q_r(0)$ and $q_r(T)$ specify the curves $(q_r(t),p_r(t))$ on $[0,T]$. However, ignoring for now $q_d(T)$, note that the initial condition $q_d(0)$ specifies $q_d(t)$ on $[0,T]$ from the initial value problem \eqref{eq:model-ham-eom-c} with $q_d(0)$ fixed, since \eqref{eq:model-ham-eom-c} is a first-order ordinary differential equation (ODE) in only the variable $q_d$. Substituting this solution curve into \eqref{eq:model-ham-eom-d}, we see that this is a time-dependent ODE in solely the variable $p_d$, but there are no specified conditions for $p_d(0)$ or $p_d(T)$. Hence, Type I boundary conditions are incomplete for this problem. Furthermore, since \eqref{eq:model-ham-eom-c} is a first-order ODE in solely the variable $q_d$, one cannot generically impose both $q_d(0)$ and $q_d(T)$; that is, with $q_d(0)$ fixed, there is only one possible choice of $q_d(T)$ that is consistent with the solution curve $q_d(t)$ given by fixing $q_d(0)$. Thus, we see that Type I boundary conditions, for the model degenerate Hamiltonian, are simultaneously underdetermined, in that they do not completely specify the solution $(q(t),p(t))$, and overdetermined, in that we cannot generically specify both $q(0)$ and $q(T)$.

\textbf{Type II/III.} In the Type II case, we have boundary conditions $q(0)= (q_{r0},q_{d0})$ and $p(T) = (p_{r1}, p_{d1})$. Again, by \Cref{assumption:complete}, $(q_r(0),p_r(T))=(q_{r0},p_{r1})$ specifies the curve $(q_r(t), p_r(t))$ on $[0,T]$. Now, the dynamics for $q_d$ \cref{eq:model-ham-eom-c} only depend on $q_d$; thus, by specifying $q_d(0) = q_{d0}$, we obtain an initial value problem involving only $q_d$ and thus, determines a curve $q_d(t)$ on $[0,T]$. Finally, by substituting this curve $q_d(t)$ into the evolution equation for $p_d$ \eqref{eq:model-ham-eom-d}, we obtain an ODE in only the variable $p_d$ with a specified value $p_d(T) = p_{d1}$; thus, this determines a curve $p_d(t)$ on $[0,T]$. The Type III case is analogous to the Type II case; simply in reverse time $s(t) := T-t$.

\textbf{Type IV.} Type IV boundary conditions are given by specifying the momenta at the initial and final time, $p(0) = p_0, p(T) = p_1$. We immediately see that these boundary conditions are incomplete for the model degenerate Hamiltonian, since \Cref{eq:model-ham-eom-c} is an ODE solely in the variable $q_d$, for which no conditions are specified. 

Summarizing the results of our discussion, for the model degenerate Hamiltonian, the completeness of the boundary conditions are:

\begin{table}[H]
\begin{center}
\begin{tabular}{c|c|c}
     Type 0 & $q(0) = q_0, p(0) = p_0$ & Complete \\ \hline
     Type I & $q(0) = q_0, q(T) = q_1$ & Incomplete \\ \hline
     Type II & $q(0) = q_0, p(T) = p_1$ & Complete \\ \hline
     Type III & $p(0) = p_0, q(T) = q_1$ & Complete  \\ \hline
     Type IV & $p(0) = p_0, p(T) = p_1$ & Incomplete
\end{tabular}
\caption{Completeness of the five types of boundary conditions for the model degenerate Hamiltonian \eqref{eq:model-hamiltonian}.}
\label{table:bcs-completeness}
\end{center}
\end{table}

\begin{remark}[Maximally Degenerate Hamiltonians] Maximally degenerate Hamiltonians are particularly interesting as they arise in the adjoint sensitivity analysis of ODEs; namely, for an ODE $\dot{q} = f(t,q)$, the adjoint Hamiltonian~\cite{TrLe2024adj} associated with this ODE is 
$$ H(t,q,p) = \langle p, f(t,q)\rangle + g(t,q), $$
which depends on a choice of a ``running cost function" $g$. This is precisely the form of a maximally degenerate Hamiltonian, i.e., adjoint Hamiltonians and maximally degenerate Hamiltonians coincide \cite{TrLe2024adj}. The associated Hamiltonian system is useful in computing the sensitivities of a cost function
$$ \mathcal{J}[q] = C(q(T)) + \int_0^T g(t,q) dt,$$
where the input $q$ of $\mathcal{J}$ is any curve satisfying the ODE $\dot{q} = f(t,q)$, $C$ is a terminal cost function and $g$ is a running cost function. 

We will discuss this application in more detail in \Cref{sec:adjoint-sensitivity}. An analogous discussion for the model Hamiltonian shows that Type 0, Type II and Type III boundary conditions are complete for maximally degenerate Hamiltonians, whereas Type I and Type IV are incomplete. 
\end{remark}

\subsection{Type II/III Variational Principles}\label{sec:type-ii-iii-vp}
From the above discussion, Type II/III boundary conditions are natural in the context of degenerate Hamiltonian systems, compared to their Type I/IV counterparts. Type II boundary conditions are of particular interest in the setting of ODE and partial differential equation (PDE) constrained optimization and optimal control, where they arise in adjoint sensitivity analysis \cite{Ca1981, CaLiPeSe2003} and Pontryagin's maximum principle \cite{BuLe2014, Bloch2015, deLe2007}. Furthermore, Type II/III variational principles provide a natural framework from which to construct Hamiltonian variational integrators \cite{LeZh2009, ScLe2017, HoTy2018, KrTy2020}; we provide an overview of such applications in \Cref{sec:generalizations}.

Motivated by this, we discuss several variational principles which incorporate Type II/III boundary conditions. Since the Type III case can be formally obtained from the Type II case from the reverse time transformation $s(t) = T-t$, we focus on the Type II case. Note that, when discussing such variational principles, we will only work with the Hamiltonian, with no reference to a Lagrangian equivalent, and make no assumptions regarding regularity of the Hamiltonian so as to include the possibility of degenerate Hamiltonians.

In \Cref{sec:type-ii-vp-vector-space}, we recall the Type II variational principle on vector spaces of \citet{LeZh2009}. As we will see, this variational principle is not intrinsic. To fix these issues, we introduce a Type II d'Alembert variational principle on vector spaces in \Cref{sec:type-ii-vp-d'Alembert}. We extend this variational principle to parallelizable manifolds in \Cref{sec:type-ii-vp-parallelizable}. Subsequently, in \Cref{sec:free-boundary-type-ii}, we show how this variational principle can be formulated globally on any manifold, at the cost of converting the Type II boundary value problem into a free boundary problem. Finally, we extend this new variational principle and the Type II variational principle of \citet{LeZh2009} to the setting of infinite-dimensional Banach spaces in \Cref{sec:inf-dim-type-ii-vp}.

\subsubsection{Type II Variational Principle on Vector Spaces}\label{sec:type-ii-vp-vector-space}
We first recall the Type II variational principle on a vector space, as introduced in \cite{LeZh2009}, where we work on the cotangent bundle $Q \times Q^* \ni (q,p)$, where again $Q$ is a finite-dimensional vector space.

To construct a variational principle $\delta S[q,p] = 0$ for Type II boundary conditions $q(0) = q_0 \in Q$, $p(T) = p_1 \in Q^*$, we define the variations as follows. For a curve $(q,p): [0,T] \rightarrow T^*Q$, consider the space of all one-parameter families of curves $(q_\epsilon, p_\epsilon): [0,T] \rightarrow Q \times Q^*$ such that
\begin{equation}
\begin{aligned}\label{eq:one-param-family-local}
    q_0(t) &= q(t), &p_0(t) &= p(t), & &\text{ for all } t \in [0,T],  \\
    q_\epsilon(0) &= q_0, &p_\epsilon(T) &= p_1, & &\text{ for all } \epsilon. 
\end{aligned}
\end{equation}
Differentiating any such one-parameter family of curves with respect to $\epsilon$ at $\epsilon = 0$, we obtain the induced variation vector field,
\begin{align*}
    \delta q(t) = \frac{d}{d\epsilon}\Big|_{\epsilon = 0} q_\epsilon(t), \\
    \delta p(t) = \frac{d}{d\epsilon}\Big|_{\epsilon = 0} p_\epsilon(t).
\end{align*}
By the first line of \eqref{eq:one-param-family-local}, $(\delta q(t),\delta p(t))$ is tangent to $(q(t),p(t))$. By the second line of \eqref{eq:one-param-family-local}, such a variation vector field satisfies $\delta q(0)=0$, $\delta p(T) = 0$. We refer to such variations as Type II variations. Note that this space of variations can be identified with the tangent space of an appropriate space of curves. Namely, for given Type II boundary conditions $q(0)=q_0, p(T) = p_1$, define the space of all $C^1$ curves satisfying these boundary conditions,
\begin{equation}\label{eq:curve-type-ii}
C^1_{q_0,p_1}([0,T], Q \times Q^*) := \{ (q,p) \in C^1([0,T],Q\times Q^*): q(0)=q_0, p(T) = p_1 \}. 
\end{equation}
Then, the space of Type II variations is the tangent space to this space of curves, given by
\begin{align}\label{eq:curve-variation-type-ii}
T&C^1_{q_0,p_1}([0,T], Q \times Q^*) \\
&= C^1_{q_0,p_1}([0,T], Q \times Q^*) \times \{ (\delta q, \delta p) \in C^1([0,T],T(Q\times Q^*)): \nonumber \\
& \qquad \qquad (\delta q(t),\delta p(t)) \in T_{(q(t),p(t))}(Q \times Q^*), \delta q(0)=0, \delta p(T) = 0 \}. \nonumber
\end{align}
Given a Hamiltonian $H:[0,T] \times Q \times Q^* \rightarrow \mathbb{R}$, define the Type II action functional $S_{\textup{II}}: C^1_{q_0,p_1}([0,T], Q \times Q^*)\rightarrow \mathbb{R}$ by
\begin{equation}\label{eq:type-ii-action-local}
    S_{\textup{II}}[q,p] := \langle p(T), q(T)\rangle - \int_0^T [\langle p(t), \dot{q}(t)\rangle - H(t,q(t),p(t))]\,dt.
\end{equation}

We define the Type II variational principle to be the stationarity of the action $\delta S_{\textup{II}}[q,p] = 0$ for all Type II variations $(\delta q, \delta p) \in T_{(q,p)} C^1_{q_0,p_1}([0,T], Q \times Q^*)$, where the variation of the action is defined as $\delta S_{\textup{II}}[q,p] = \frac{d}{d\epsilon}\big|_{\epsilon = 0} S_{\textup{II}}[q_\epsilon,p_\epsilon]$. 

\begin{prop}[Type II Variational Principle on a Vector Space]\label{prop:type-ii-vp-local}
    The Type II action is stationary $\delta S_{\textup{II}} = 0$ for all Type II variations if and only if Hamilton's equations with Type II boundary conditions hold:
    \begin{subequations}\label{eq:type-ii-hamilton-equations}
    \begin{alignat}{2}
        \dot{q}(t) &= D_pH(t,q,p),&\, &\text{ for all } t\in(0,T), \label{eq:type-ii-hamilton-equations-a}\\
        \dot{p}(t) &= -D_qH(t,q,p),&\, &\text{ for all } t\in(0,T), \label{eq:type-ii-hamilton-equations-b}\\
        q(0) &= q_0, \label{eq:type-ii-hamilton-equations-c}\\
        p(T) &= p_1. \label{eq:type-ii-hamilton-equations-d}
    \end{alignat}
    \end{subequations}
    \begin{proof}
        The ``if" direction is trivial to verify. For the ``only if" direction, note that the space of curves \eqref{eq:curve-type-ii} enforces the boundary conditions \eqref{eq:type-ii-hamilton-equations-c}-\eqref{eq:type-ii-hamilton-equations-d} and these are preserved by Type II variations \eqref{eq:curve-variation-type-ii} by construction. So, we simply verify that the stationarity condition yields Hamilton's equations. We explicitly compute the variation of the action to be
        \begin{align*}
            \delta S_{\textup{II}}[q,p]  &= \langle \cancel{\delta p(T)}, q(T)\rangle + \langle p(T), \delta q(T)\rangle \\
            &\qquad - \int_0^T \left[\langle \delta p(t), \dot{q}(t) - D_p H \rangle + \langle p(t), \frac{d}{dt}\delta q(t)\rangle - \langle D_qH, \delta q\rangle \right]\,dt \\
            &= \langle p(T),\delta q(T)\rangle - \langle p(T), \delta q(T)\rangle + \langle p(0), \cancel{\delta q(0)}\rangle \\
            & \qquad - \int _0^T \left[\langle \delta p(t), \dot{q}(t) - D_p H \rangle + \langle -\dot{p} - D_qH, \delta q \rangle \right]\,dt.
        \end{align*}
        A standard application of the fundamental lemma of the calculus of variations yields Hamilton's equations \eqref{eq:type-ii-hamilton-equations-a}-\eqref{eq:type-ii-hamilton-equations-b} for $t \in (0,T)$.
    \end{proof}
\end{prop}

\Cref{prop:type-ii-vp-local} provides a Type II variational principle when the configuration space is a vector space. In moving to the more general case of a manifold $M$, there are two fundamental issues with the above argument. First, in the definition of the Type II action \eqref{eq:type-ii-action-local}, we see a pairing of the form $\langle p(T), q(T)\rangle$. While this pairing is well-defined in the vector space case since $p(T) \in Q^*$, $q(T) \in Q$, this will not generally make sense in the manifold case, since $p(T) \in T^*_{q(T)}M$ and $q(T) \in M$, i.e., this would be a pairing between a cotangent vector and a base point. This leads to the second issue: the Type II boundary condition is given by specifying $q(0) = q_0 \in M$ and $p(t) = p_1 \in T^*_{q(T)}M$, yet the base point $q(T)$ is unspecified; one cannot in general specify a cotangent vector $p_1$ without specifying its base point. Thus, the space of curves \eqref{eq:curve-type-ii}, with $Q \times Q^*$ replaced by $T^*M$, would generally not make sense. This was not an issue in the vector space case, since the cotangent bundle $T^*Q$ is globally trivializable as $Q \times Q^*$ so we can specify the cotangent vector $p_1$ in $Q^*$. 

We will start by resolving the first issue in the vector space setting, and then discuss how the second issue can be resolved when $M$ is parallelizable. Subsequently, we show that the second issue can be resolved globally if we allow the Type II Hamiltonian system to be a free boundary problem. 

\subsubsection{Type II d'Alembert Variational Principle on Vector Spaces} \label{sec:type-ii-vp-d'Alembert}
\begin{figure}[H]
\begin{center}
\includegraphics[width=120mm]{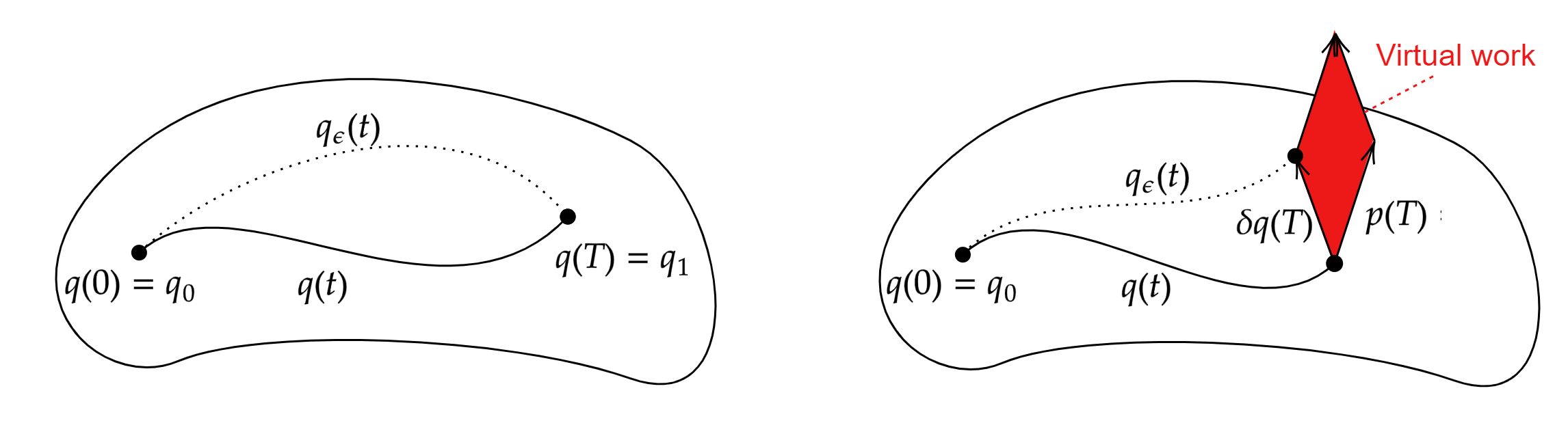}
\caption{Schematic of Type I variations (left) and partial variations (right) on the base manifold.} 
\label{figure:type-ii-virtual-work}
\end{center}
\end{figure}

We remain in the vector space case as before but modify the action and the corresponding variational principle. Since we do not want the pairing $\langle p(T), q(T)\rangle$ appearing in the action because it represents a pairing between a cotangent vector $p(T) \in Q^*$ and a position $q(T) \in Q$, we will remove this term which returns us to the classical Hamiltonian action \eqref{eq:classical-hamiltonian-action},
$$ S[q,p] = \int_0^T [\langle p(t), \dot{q}(t)\rangle - H(t,q,p)]\,dt. $$
We would like to define an appropriate space of variations. To do so, recall the variation of the above action is, for arbitrary $(\delta q,\delta p)$,
\begin{align*}
        \delta S[q,p] &=   \langle p(T), \delta q(T)\rangle - \langle p(0), \delta q(0)\rangle  \\
        & \qquad + \int _0^T \left[\langle \delta p(t), \dot{q}(t) - D_p H \rangle + \langle -\dot{p} - D_qH, \delta q \rangle \right]\,dt.
\end{align*}
In order to obtain Hamilton's equations of motion from the integral terms, the boundary terms $\langle p(0), \delta q(0)\rangle$ and $\langle p(T), \delta q(T)\rangle$ have to be made to vanish. For the former term, since we enforce fixed initial position $q(0)=q_0$, we restrict the space of variations to satisfy $\delta q(0) = 0$, as before. For the latter term, we do not have the freedom to choose $\delta q(T)$ since we do not fix the terminal position $q(T)$. Thus, since $\delta q(T)$ is arbitrary, stationarity can only be obtained with $p(T) = 0$, since otherwise $\langle p(T), \delta q(T)\rangle$ can be made arbitrarily large. We see the emergence of a natural boundary condition for $p(T)$, i.e., one that arises from the variational principle, as opposed to the essential boundary condition for $q(0)$ which arises from restricting the variations \textit{a priori} (the terminology ``essential" and ``natural" boundary conditions are standard terminology in weak formulations of boundary value problems, see, e.g., \cite{BrSc2008, VuAt2008}). We thus define the space of partial variations as variations $(\delta q, \delta p)$ satisfying $\delta q(0)=0$. A schematic of the Type I variations fixing both endpoints $q(0)=q_0, q(T) = q_1$ and partial variations fixing only $q(0)=q_0$ on the base manifold is shown in \Cref{figure:type-ii-virtual-work}.
Since Type I variations satisfy $\delta q(0)=0=\delta q(T)$, the boundary terms introduce no ``virtual work", $\langle p(T),\delta q(T)\rangle = 0 = \langle p(0), \delta q(0)\rangle$, and thus, Hamilton's equations can be obtained from the stationarity condition $\delta S = 0$ of the classical action \eqref{eq:classical-hamiltonian-action} as discussed in \Cref{sec:mechanics}. On the other hand, for partial variations, $\delta q(T)$ is arbitrary so virtual work can be introduced through variations (represented by the area of the shaded region in \Cref{figure:type-ii-virtual-work}); as such, for the classical action \eqref{eq:classical-hamiltonian-action}, stationarity of the action with respect to all partial variations holds only if the natural boundary condition $p(T) = 0$ holds.

To state the variational principle precisely, for a given $q_0 \in Q$, we define the space of curves fixing $q(0)=q_0$ as
\begin{equation}\label{eq:curve-partial}
C^1_{q_0}([0,T], Q \times Q^*) := \{ (q,p) \in C^1([0,T],Q\times Q^*): q(0)=q_0 \};
\end{equation}
then, the space of partial variations is given by the tangent space to this space of curves
\begin{align}\label{eq:curve-variation-partial}
T C^1_{q_0}&([0,T], Q \times Q^*) \\ 
&= C^1_{q_0}([0,T], Q \times Q^*)  \times \{ (\delta q, \delta p) \in C^1([0,T],T(Q\times Q^*)): \nonumber \\
&\qquad \qquad (\delta q(t),\delta p(t)) \in T_{q(t),p(t)}(Q\times Q^*),\delta q(0)=0 \}. \nonumber
\end{align}
Furthermore, to allow for an inhomogeneous natural boundary condition $p(T) = p_1 \in Q^*$, we must modify the variational principle. As noted before, if $p(T) = p_1 \neq 0$, the variation of the action can be made arbitrarily large through the virtual work $\langle p(T), \delta q(T)\rangle$. Thus, there is no hope for stationarity of $S$ if we demand $p(T) = p_1 \neq 0$. Instead, we will demand that the variation of the action $\delta S[q,p]$ under partial variations is equal to the virtual work done by $p_1$ under the variation $\delta q(T)$, i.e., $\langle p_1,\delta q(T)\rangle$. 

For fixed $q_0 \in Q$, fixed $p_1 \in Q^*$ and defining the action as a mapping $S: C^1_{q_0}([0,T], Q \times Q^*) \rightarrow \mathbb{R}$, we define the Type II d'Alembert variational principle to be the condition
\begin{equation}\label{eq:d'Alembert-variational-principle}
    \delta S[q,p] = \langle p_1, \delta q(T)\rangle,
\end{equation}
for all partial variations $(\delta q,\delta p) \in T_{(q,p)} C^1_{q_0}([0,T], Q \times Q^*)$.

\begin{remark}
    We refer to this as a d'Alembert variational principle due to its similarity in form with the Lagrange--d'Alembert variational principle for forced Lagrangian systems (for a summary of this principle, see \cite{VuAt2008}; in the context of variational integrators, see \cite{LeMaOrWe2004}). 
\end{remark}

\begin{prop}[Type II d'Alembert Variational Principle]\label{prop:type-ii-vp-d'Alembert}
    The above Type II d'Alembert variational principle holds if and only if Hamilton's equations with Type II boundary conditions hold:
    \begin{equation}\label{eq:type-ii-hamilton-equation-d'Alembert}
    \begin{aligned}
        \dot{q}(t) &= D_pH(t,q,p), & \dot{p}(t) &= -D_qH(t,q,p), & \text{ for all } t \in (0,T),\\
        q(0) &= q_0, & p(T) &= p_1. 
    \end{aligned}
    \end{equation}
    \begin{proof}
        An analogous computation to the homogeneous case considered before shows that the Type II d'Alembert variational principle is equivalent to \eqref{eq:type-ii-hamilton-equation-d'Alembert}. The only difference in the calculation is that now the problematic boundary term takes the form $\langle p(T) - p_1, \delta q(T)\rangle$, which vanishes for all $\delta q(T)$ if and only if   $p(T) = p_1$; thus, the terminal momentum condition $p(T)=p_1$ again arises as a natural boundary condition. 
    \end{proof}
\end{prop}

\begin{remark}\label{remark:intrinsic-type-ii-vector-space}
We can state the Type II d'Alembert variational principle intrinsically as follows. Let $\Theta$ be the canonical one-form on $Q \times Q^*$ with coordinate expression $\Theta = p\, dq$. Let $\pi_{Q \times Q^*}: [0,T] \times Q \times Q^* \rightarrow Q \times Q^*$ be the standard projection onto the second and third factors. Define the Cartan form associated with a Hamiltonian $H: [0,T] \times Q \times Q^* \rightarrow \mathbb{R}$ as
$$ \Theta_H := \pi_{Q \times Q^*}^*\Theta - H dt, $$
which is a one-form on $[0,T] \times Q \times Q^*$. Consider the space of curves on $[0,T] \times Q \times Q^*$ covering the identity on $[0,T]$, i.e., curves of the form $\psi(t) = (t,q(t),p(t))$. In terms of these quantities, the classical action \eqref{eq:classical-hamiltonian-action} can be expressed as
\begin{equation}\label{eq:hamiltonian-cartan-action}
    S[\psi] = \int_0^T \psi^* \Theta_H.
\end{equation}
Let $\pi_{[0,T]}: [0,T] \times Q \times Q^* \rightarrow [0,T]$, $\pi_{Q}: [0,T] \times Q \times Q^* \rightarrow Q$, and $\pi_{Q^*}: [0,T] \times Q \times Q^* \rightarrow Q^*$ be the standard projections onto the first, second, and third factors, respectively. For a given curve $\psi$, we define the partial variations as the space of time-dependent vector fields $X$ tangent to $\psi$ which are vertical with respect to $\pi_{[0,T]}$ and satisfy $T\pi_Q X(0,\psi(0)) = 0$; note that this gives an intrinsic definition of the space of partial variations \eqref{eq:curve-variation-partial}. Let $\psi_1$ be the curve with the same position as $\psi$ but fixed momentum coordinate $p_1$, i.e., $\psi_1(t) = (t,q(t),p_1)$. The Type II d'Alembert variational principle \eqref{eq:d'Alembert-variational-principle} can then equivalently be expressed as 
\begin{equation}\label{eq:d'Alembert-variational-principle-cartan}
    d S[\psi]\cdot X = \psi_1^*(i_X \Theta_H)(T),
\end{equation}
for all partial variations $X$. Note that the coordinate expression for the right hand side of \eqref{eq:d'Alembert-variational-principle-cartan} with $X = \delta q \frac{\partial}{\partial q} + \delta p \frac{\partial}{\partial p}$ is precisely the right hand side of \eqref{eq:d'Alembert-variational-principle}, i.e., is an intrinsic expression for the virtual work done by $p_1$ through the variation $\delta q(T)$. Finally, Hamilton's equations follow from a direct calculation. Let $\Omega_H = -d\Theta_H$ and $\varphi_\epsilon$ be the time-$\epsilon$ flow of $X$; then,  we have
\begin{align*}
    \psi^*_1(i_X \Theta_H)(T) &= d S[\psi]\cdot X = \frac{d}{d\epsilon}\Big|_{\epsilon = 0} S[\varphi_\epsilon \circ \psi] = \int_0^T \psi^* \frac{d}{d\epsilon}\Big|_{\epsilon = 0} \varphi_{\epsilon}^* \Theta_H \\
    &= \int_0^T \psi^* \mathcal{L}_X \Theta_H = - \int_0^T \psi^*(i_X\Omega_H) + \int_0^T \psi^* d(i_X\Theta_H) \\
    &= - \int_0^T \psi^*(i_X\Omega_H) + \int_0^T d(\psi^*( i_X\Theta_H)) \\
    &= - \int_0^T \psi^*(i_X\Omega_H) + \psi^* (i_X \Theta_H)(T) - \psi^*(i_X\Theta_H)(0).
\end{align*}
The last term vanishes since $\psi^*(i_X\Theta_H)(0) = \langle \pi_{Q^*}\psi(0), T\pi_QX(\psi(0))\rangle = 0$ and thus, comparing the two remaining boundary terms in the above equation, we again obtain natural boundary conditions $p(T)=p_1$, so that the above reduces to $\int_0^T \psi^*(i_X\Omega_H) = 0$. A standard application of the fundamental lemma of the calculus of variations yields the implicit Hamilton's equations $\psi^*(i_X\Omega_H) = 0$ for all partial variations $X$.
\end{remark}

We have thus remedied the first issue with the traditional Type II variational principle (\Cref{prop:type-ii-vp-local}) containing a term of the form $\langle p(T), q(T)\rangle$. However, the second issue still remains; namely, that specifying the cotangent vector $p_1$ without specifying its base point $q_1$ does not make intrinsic sense on a manifold. In essence, we will need additional structure to make intrinsic sense of specifying the cotangent vector $p_1$. We will now discuss how to formulate this variational principle when the manifold is parallelizable and subsequently, we show how to do this globally without this assumption by relaxing the condition that the cotangent vector $p_1$ is fixed and instead considering the Type II Hamiltonian system as a free boundary problem. 

\subsubsection{Type II Variational Principle for Parallelizable Manifolds}\label{sec:type-ii-vp-parallelizable}
Consider the case where $M$ is parallelizable, i.e., the tangent bundle $\pi_{TM}:TM\rightarrow M$ is trivializable. That is, we assume there exists a trivialization $TM \cong M \times V$, where $V$ is a vector space, given by a trivialization map $\Phi$ which is a diffeomorphism,
\begin{align*}
    \Phi:  M \times V \rightarrow TM,
\end{align*}
which is fiber-preserving $\Phi(q,\cdot): V \rightarrow T_qM$ for all $q \in M$ and furthermore, $\Phi(q,\cdot): V \rightarrow T_qM$ is a linear isomorphism for all $q \in M$. We write the action of this isomorphism on $v \in V$ as $\Phi(q)v := \Phi(q,v)$. Similarly, we denote the inverse isomorphism by $\Phi(q)^{-1}: T_qM \rightarrow V $. 
By duality, the trivialization of $TM$ induces a trivialization of the cotangent bundle $T^*M$, where now the isomorphisms are given by
\begin{align*}
    \Phi(q)^{-*} := (\Phi(q)^{-1})^*&: V^* \rightarrow T_q^*M, \\
    \Phi(q)^*&: T_q^*M \rightarrow V^*. 
\end{align*}
The derivative of $\Phi(q)$ and its inverse will be useful in our discussion below. Letting $L(A,B)$ denote the space of linear maps between vector spaces $A$ and $B$, we denote the derivative as
\begin{equation*}
\begin{aligned}
    D\Phi(q): T_qM &\rightarrow L(V, T_qM)  \\
        v &\mapsto D\Phi(q)\cdot v, \\
    D\Phi(q)\cdot v: V &\rightarrow T_qM \\
                    \xi &\mapsto D\Phi(q)\cdot v \cdot \xi,
\end{aligned}
\end{equation*}
and similarly for its inverse $D(\Phi(q)^{-1}): T_qM \rightarrow L(T_qM,V)$. 

\begin{remark}\label{rmk:derivative-map}
    More precisely, the derivative above is the differential of a smooth map between manifolds. For concreteness, let $U$ be a local chart about $q \in M$ for which we have the trivialization $TU \cong U \times V$. Let $\{q^i\}$ be local coordinates on $U$ and let $\{e_j\}$ be a basis for $V$; then, the isomorphism $\Phi(q): V \rightarrow T_qM$ can be represented locally as a matrix $\Phi(q^i)^j_k$, whose components depend on $q = (q^i)$, relative to this basis of $V$. The derivative map $D\Phi(q)$ then has the local expression 
    $$ \frac{\partial}{\partial q^l} \Phi(q^i)^j_k, $$
    and similarly $D\Phi(q)\cdot v$ has the local expression
    $$ v^l\frac{\partial}{\partial q^l} \Phi(q^i)^j_k, $$
    where $v \in T_qM$ has the expression $v^l e_l \in V$ relative to the trivialization $TU \cong U \times V$.
\end{remark}

Now, we would like to enforce Type II boundary conditions for a Hamiltonian system. Using the trivialization of $T^*M$, we can supply trivialized Type II boundary conditions of the form $q(0) = q_0 \in M$ and $\mu(T) = \mu_1 \in V^*$ where we regard $(q,\mu) \in M \times V^*$ as the trivialization of $(q,p) = (q,\Phi(q)^{-*}\mu) \in T^*_q M$. Associated with a Hamiltonian $H:[0,T] \times T^*M \rightarrow \mathbb{R}$, we define the trivialized Hamiltonian
\begin{align}\label{eq:trivialized-hamiltonian}
    h &:\, [0,T] \times M \times V^* \rightarrow \mathbb{R} \\
                h(t,q,\mu) &= H(t,q,\Phi(q)^{-*}\mu). \nonumber
\end{align}
Simlarly, we define the trivialized action corresponding to \eqref{eq:classical-hamiltonian-action}, noting that $\langle p,\dot{q}\rangle = \langle \Phi(q)^{-*}\mu, \dot{q}\rangle = \langle \mu, \Phi(q)^{-1} \dot{q}\rangle$, where the duality pairing on the right hand side is the duality pairing on $V^* \times V$. Thus, the trivialized action is
\begin{align}\label{eq:trivialized-action}
    s[q,\mu] := \int_0^T [\langle \mu, \Phi(q)^{-1} \dot{q}\rangle - h(t,q,\mu) ]\, dt.
\end{align}
We will again utilize the d'Alembert variational principle, noting that the Type II action functional \eqref{eq:type-ii-action-local} does not make intrinsic sense due to the pairing $\langle p(T), q(T)\rangle$ since $M$ may not generally be a vector space. Corresponding to trivialized Type II boundary conditions, we require the trivialized variation $\eta := \Phi(q)^{-1}\delta q$ to satisfy $\eta(0) = 0$ and allow the variation $\delta \mu$ to be arbitrary.

More precisely, for a given $q_0 \in M$, we define the space of curves fixing $q(0)=q_0$ as
\begin{equation}\label{eq:curve-partial-trivialized}
C^1_{q_0}([0,T], M \times V^*) := \{ (q,\mu) \in C^1([0,T],M\times V^*): q(0)=q_0 \};
\end{equation}
then, by trivialization, the space of trivialized partial variations is seen to be
\begin{equation}\label{eq:curve-variation-partial-trivialized}
T C^1_{q_0}([0,T], M \times V^*) = C^1_{q_0}([0,T], M \times V^*)  \times \{ (\eta,\delta \mu) \in C^1_{q_0}([0,T], V \times V^*), \eta(0) = 0 \}.
\end{equation}
As in the previous sections, we interpret the action as a functional on this space of curves, $s:C^1_{q_0}([0,T], M \times V^*) \rightarrow \mathbb{R}$.

This gives us the following Type II d'Alembert variational principle for a parallelizable manifold.

\begin{prop}
    The Type II d'Alembert variational principle,
    $$ \delta s[q,\mu] = \langle \mu_1(T), \eta(T)\rangle, $$
    for all trivialized partial variations $(\eta, \delta \mu) \in T_{(q,\mu)} C^1_{q_0}([0,T], M \times V^*)$ holds if and only if the Hamel equations on $M \times V^*$ with trivialized Type II boundary conditions hold:
    \begin{subequations}\label{eq:type-ii-hamel-equations}
    \begin{alignat}{2}
        \dot{q}(t) &= \Phi(q) D_\mu h(t,q,\mu),&\, &\text{ for all } t\in(0,T), \label{eq:type-ii-hamel-equations-a} \\
        \dot{\mu}(t) &= \textup{ad}^*_\xi \mu - \Phi(q)^* D_qh(t,q,\mu),&\, &\text{ for all } t\in(0,T), \label{eq:type-ii-hamel-equations-b} \\
        q(0) &= q_0, \label{eq:type-ii-hamel-equations-c} \\
        \mu(T) &= \mu_1, \label{eq:type-ii-hamel-equations-d}
    \end{alignat}
    where $\xi := \Phi(q)^{-1}\dot{q}$ denotes the trivialized velocity and $\textup{ad}^*_\xi$ will be defined in the proof. 
    \end{subequations}
    \begin{proof}
        We compute the variation of the action
        \begin{align*}
             \delta s[q,\mu] &= \int_0^T \Big[\langle \delta \mu, \Phi(q)^{-1} \dot{q} - D_\mu h(t,q,\mu)\rangle + \langle \mu, \delta (\Phi(q)^{-1} \dot{q}) \rangle - \langle D_qh(t,q,\mu),\delta q\rangle \Big]\, dt \\
            &= \int_0^T \langle \delta \mu, \Phi(q)^{-1} \dot{q} - D_\mu h(t,q,\mu)\rangle\,dt \\
                &\qquad + \int_0^T \Big[\langle \mu, D(\Phi(q)^{-1})\cdot \delta q \cdot \dot{q} + \Phi(q)^{-1} \frac{d}{dt} \delta q\rangle - \langle D_qh(t,q,\mu), \Phi(q) \eta \rangle \Big]\,dt \\
            &= \int_0^T \langle \delta \mu, \Phi(q)^{-1} \dot{q} - D_\mu h(t,q,\mu)\rangle\,dt \\
                &\qquad + \int_0^T \Big[\langle \mu, \underbrace{D(\Phi(q)^{-1})\cdot \delta q \cdot \dot{q} + \Phi(q)^{-1} \frac{d}{dt} \delta q\rangle}_{=:\ (\textbf{a})} - \langle \Phi(q)^* D_qh(t,q,\mu), \eta \rangle \Big]\,dt.
        \end{align*}
        We express $(\textbf{a})$ in terms of the trivialized variation $\eta$ and the trivialized velocity $\xi$ as follows,
        \begin{align*}
            (\textbf{a}) & = D(\Phi(q)^{-1})\cdot \delta q \cdot \dot{q} + \Phi(q)^{-1} \frac{d}{dt} \Phi(q) \eta \\
            &=  D(\Phi(q)^{-1})\cdot \delta q \cdot \dot{q} + \Phi(q)^{-1} D\Phi(q) \cdot \dot{q} \cdot \eta + \frac{d\eta}{dt} \\
            &=  - \Phi(q)^{-1} D\Phi(q)\cdot \delta q \cdot \Phi^{-1}\dot{q} + \Phi(q)^{-1} D\Phi(q) \cdot \dot{q} \cdot \eta + \frac{d\eta}{dt} \\
            &= - \Phi(q)^{-1} D\Phi(q)\cdot \Phi(q) \eta \cdot \xi + \Phi(q)^{-1} D\Phi(q) \cdot \Phi(q) \xi \cdot \eta + \frac{d\eta}{dt},
        \end{align*}
        where we used the identity for the derivative of the inverse (recalling our notation for the derivative of $\Phi^{-1}$),
        $$D(\Phi(q)^{-1})\cdot u \cdot v = - \Phi(q)^{-1} D\Phi(q) \cdot u \cdot \Phi(q)^{-1} v.$$
        Define the antisymmetric bracket, for $u,v \in V$,
        \begin{equation}\label{eq:hamel-bracket}
            [u,v]_{q} :=  \Phi(q)^{-1} D\Phi(q) \cdot \Phi(q) u \cdot v - \Phi(q)^{-1} D\Phi(q)\cdot \Phi(q) v \cdot u.
        \end{equation}
        As observed in \cite{Bloch2015, Bloch2009a, Bloch2009b}, \eqref{eq:hamel-bracket} defines a Lie bracket on $V$, which depends on $q \in Q$, induced from the Lie bracket of vector fields on $TM$ (see \Cref{remark:hamel}).
        
        Then, $(\textbf{a})$ can be expressed as
        $$ (\textbf{a}) = \frac{d\eta}{dt} + [\xi,\eta]_q. $$
        By duality, we define the coadjoint action $\textup{ad}^*_{\xi}$ to be the adjoint of $[\xi,\cdot]_q: V \rightarrow V$, i.e., for $\alpha \in V^*$,
        $$\langle \textup{ad}^*_\xi \alpha, v \rangle := \langle \alpha, [\xi, v]_q \rangle,  \text{ for all } v \in V. $$

        Returning to the variation of the action, we have from integrating by parts,
        \begin{align*}
            \delta s[q,\mu] &= \int_0^T \langle \delta \mu, \Phi(q)^{-1} \dot{q} - D_\mu h(t,q,\mu)\rangle\,dt \\
            & \qquad + \int_0^T \Big[\left\langle \mu , \dot{\eta} + [\xi,\eta]_q \right\rangle - \langle \Phi(q)^* D_qh(t,q,\mu), \eta \rangle \Big]\,dt \\
            &= \langle \mu(T), \eta(T) \rangle - \langle \mu(0), \eta(0)\rangle + \int_0^T \langle \delta \mu, \Phi(q)^{-1} \dot{q} - D_\mu h(t,q,\mu)\rangle\,dt \\
                &\qquad + \int_0^T \left\langle -\dot{\mu} + \textup{ad}^*_\xi \mu - \Phi(q)^* D_qh(t,q,\mu),\eta \right\rangle dt.
        \end{align*}
        The Type II d'Alembert variational principle $\delta s[q,\mu] = \langle \mu_1(T), \eta(T)\rangle $ for any variation $\eta$ with $\eta(0) = 0$ and arbitrary $\delta \mu$ is thus
        \begin{align*}
            \langle \mu(T) - \mu_1, \eta(T) \rangle &+ \int_0^T \langle \delta \mu, \Phi(q)^{-1} \dot{q} - D_\mu h(t,q,\mu)\rangle\,dt \\
            & \qquad + \int_0^T \left\langle -\dot{\mu} + \textup{ad}^*_\xi \mu - \Phi(q)^* D_qh(t,q,\mu),\eta \right\rangle dt = 0.
        \end{align*} 
        From this expression, it immediately follows that the variational principle holds if and only if Hamel's equations with Type II boundary conditions hold.
    \end{proof}
\end{prop}

    \begin{remark}\label{remark:hamel}
        Hamel's equation, see \cite{Bloch2015, Bloch2009a, Bloch2009b}, in the context of Lagrangian mechanics, arise as the transformations of the Euler--Lagrange equations under coordinate transformations mapping velocities to quasi-velocities, i.e., non-natural coordinates on $TQ$. Locally, our situation is analogous, since the trivialization of $TM$ and $T^*M$ can be thought of locally as coordinate transformations from velocity and momentum coordinates to quasi-velocity and non-canonical coordinates, respectively. The only difference is that we are working in the Hamiltonian picture with no assumption of a Lagrangian equivalent. 

        As previously mentioned, the bracket \eqref{eq:hamel-bracket} defines a Lie bracket on $V$ induced from the Lie bracket of vector fields on $TM$, as shown in \cite{Bloch2015, Bloch2009a, Bloch2009b}. In \cite{Bloch2015, Bloch2009a, Bloch2009b}, they work in coordinates, particularly with velocity coordinate transformations of the form
        $$ u_i(q) = \psi^j_i(q) \frac{\partial}{\partial q^j}, $$
        where $\{u_i\}$ is the quasi-velocity basis and $\{\partial/\partial q^j\}$ is the natural coordinate basis on $TM$, and show that the Lie bracket has the form
        \begin{align*}
            [u_i,u_j] &= c^k_{ij}(q)u_k, \\
            c^k_{ij} &= (\psi^{-1})^k_m \left[ \frac{\partial \psi^m_j}{\partial q^l} \psi^l_i - \frac{\partial \psi^m_i}{\partial q^l} \psi^l_j \right].
        \end{align*}
        This is precisely the coordinate expression for \eqref{eq:hamel-bracket}, noting that the local coordinate expression for $\Phi(q)^{-1}$ is $\psi^j_i$ (see \Cref{rmk:derivative-map}). We prefer to work in index-free notation, but refer the reader to these references for a more detailed discussion involving such coordinate transformations. 
    \end{remark}

    \begin{example}
        As an example, consider the case when $M = G$ is a Lie group. Lie groups are trivializable by left-translation, $TG \cong G \times \mathfrak{g}$ and $T^*G \cong G \times \mathfrak{g}^*$, where $\mathfrak{g} = T_eG$ denotes the Lie algebra of $G$ and $\mathfrak{g}^*$ its dual. Letting $L_g: G \rightarrow G$ denote left-translation by $g \in G$, $L_g(x) = gx$, the inverse of the trivialization maps are given by
        \begin{align*}
            \Phi(g)^{-1} &= T_g L_{g^{-1}} : T_gG \rightarrow T_eG = \mathfrak{g}, \\
            \Phi(g)^* &= T^*_eL_g: T^*_g G \rightarrow T^*_eG = \mathfrak{g}^*.
        \end{align*}
        In this case, Hamel's equations \eqref{eq:type-ii-hamel-equations-a}-\eqref{eq:type-ii-hamel-equations-b} become the Euler--Arnold equations (see, for example, \cite{deDiego2018, CuBa1997}), where the antisymmetric bracket introduced above is the Lie bracket on $\mathfrak{g}$.
    \end{example}

    We have defined a Type II d'Alembert variational principle globally on parallelizable manifolds. On the trivializing space $M \times V^*$, the Type II d'Alembert variational principle yields Hamel's equations with trivialized Type II boundary conditions $(q_0, \mu_1)$. Note that the above argument also applies locally to any (not necessarily parallelizable) manifold $M$, since its tangent and cotangent bundles are locally trivializable. Thus, we have also defined a local Type II d'Alembert variational principle for an arbitrary manifold and arbitrary coordinate chart, whereas the previous Type II variational principles on vector spaces discussed in \Cref{sec:type-ii-vp-vector-space} and \Cref{sec:type-ii-vp-d'Alembert} use canonical coordinates. Now, we show how the Type II variational principle can be formulated globally if we relax the condition that $p_1$ is a fixed cotangent vector.

    \subsubsection{Global Type II Variational Principle as a Free Boundary Problem}\label{sec:free-boundary-type-ii}
    As previously mentioned, the obstruction to formulating a Type II variational principle globally on $T^*M$ for a generic manifold $M$ is that we cannot specify a particular cotangent vector $p_1$ without specifying its basepoint $q(T)$, which is variable under the dynamics of the Hamiltonian system. We will relax the requirement of specifying a particular cotangent vector and instead let $p_1$ be variable. More precisely, we will consider $p_1$ as assigning to each $q \in M$ a cotangent vector $p_1(q) \in T^*_qM$, i.e., $p_1: M \rightarrow T^*M$ is now a section of the cotangent bundle. 
    
    Analogous to the previous sections, for $q_0 \in M$, we define the space of curves fixing $q(0)=q_0$ as
    \begin{equation}\label{eq:curve-partial-manifold}
    C^1_{q_0}([0,T], T^*M) := \{ z \in C^1([0,T], T^*M):  \pi_{T^*M}(z(0)) = q_0 \},
    \end{equation}
    where $\pi_{T^*M}: T^*M \rightarrow M$ is the cotangent bundle projection. Then, the space of partial variations is given by the tangent space to this space of curves
    \begin{align}\label{eq:curve-variation-partial-manifold}
    T C^1_{q_0}&([0,T], T^*M) \\
    & = \{ X \in C^1([0,T],TT^*M): \pi_{TT^*M}(X(0))=q_0, T\pi_{T^*M}(X(0))=0 \}, \nonumber
    \end{align}
    where $\pi_{TT^*M}: TT^*M \rightarrow M$ is the bundle projection on $TT^*M$ as a bundle over $M$. In bundle coordinates $X = (q,p,\delta q, \delta p)$ on $TT^*M$, the conditions $\pi_{TT^*M}(X(0))=q_0$ and $T\pi_{T^*M}(X(0))=0$ are simply $q(0)=q_0$ and $\delta q(0)=0$, respectively.
    
    First, we will state the variational principle in coordinates and subsequently, state it intrinsically in \Cref{remark:intrinsic-manifold-vp}. For a given section $p_1$ of $T^*M$, we define the free boundary Type II d'Alembert variational principle to be
    \begin{equation}\label{eq:free-type-ii-d'Alembert-vp}
        \delta S[q,p] = \langle p_1(q(T)), \delta q(T) \rangle,
    \end{equation}
    for all partial variations $(\delta q,\delta p)$ satisfying $\delta q(0)=0$, where the action is a functional on the space of curves $S: C^1_{q_0}([0,T], T^*M) \rightarrow \mathbb{R}$.

    A similar calculation as before yields the following result.
    \begin{prop}
    The above variational principle yields the system, in coordinates,
    \begin{subequations}\label{eq:free-type-ii-hamilton-equations}
    \begin{alignat}{2}
        \dot{q}(t) &= D_pH(t,q,p),&\,& \text{ for all } t\in(0,T), \label{eq:free-type-ii-hamilton-equations-a}\\
        \dot{p}(t) &= -D_qH(t,q,p),&\,& \text{ for all } t\in(0,T), \label{eq:free-type-ii-hamilton-equations-b}\\
        q(0) &= q_0, \label{eq:free-type-ii-hamilton-equations-c}\\
        p(T) &= p_1(q(T)). \label{eq:free-type-ii-hamilton-equations-d}
    \end{alignat}
    \end{subequations}
    \end{prop}

    \begin{remark}\label{remark:intrinsic-manifold-vp}
        \sloppy Using similar notation to \Cref{remark:intrinsic-type-ii-vector-space}, this variational principle can also be stated intrinsically as
        $$ dS[\psi]\cdot X = i_X \Theta_H(T,p_1(\pi_{M}(\psi))), $$
        for all partial variations $X$, where $\pi_M: [0,T] \times T^*M \rightarrow M$ is the natural time-dependent bundle projection. A similar calculation to \Cref{remark:intrinsic-type-ii-vector-space} shows that this variational principle is equivalent to the implicit Hamilton's equations $\psi^*(i_X\Omega_H) = 0$ with fixed initial condition $q(0)=q_0$ and an implicit boundary condition $\psi^*(i_X\Theta_H)(T) = i_X \Theta_H(T,p_1(\pi_{M}(\psi)))$, for all partial variations $X$. The coordinate expression for this implicit boundary condition is $\langle p(T), \delta q(T)\rangle = \langle p_1(q(T)), \delta q(T)\rangle$ which is the implicit form of \eqref{eq:free-type-ii-hamilton-equations-d}.
    \end{remark}

    We interpret the system \eqref{eq:free-type-ii-hamilton-equations-a}-\eqref{eq:free-type-ii-hamilton-equations-d} as a free boundary problem, since the boundary point $q(T)$ is not fixed but rather coupled to the dynamics of the system. The analysis of free boundary problems is in general complicated and problem dependent (see, for example, \cite{FreeBnd1, FreeBnd2, FreeBnd3}) so we will not make any generic statements regarding the existence and uniqueness of solutions to this system. 
    
    However, we note that for maximally degenerate Hamiltonians in particular, there always exists a unique solution to the above system. That is, for a time-dependent vector field $f: [0,T] \times M \rightarrow TM$ and a function $g: [0,T] \times M \rightarrow \mathbb{R}$, consider a Hamiltonian $H_g: [0,T] \times T^*M \rightarrow \mathbb{R}$ of the form
    $$ H_g(t,q,p) = \langle p, f(t,q)\rangle + g(t,q),$$
    (for an intrinsic definition of $H_g$, see \cite{TrLe2024adj}). Then, the free boundary problem \eqref{eq:free-type-ii-hamilton-equations} for this Hamiltonian becomes
    \begin{align*}
        \dot{q}(t) &= f(t,q), \\
        \dot{p}(t) &= -[D_qf(t,q)]^* p - D_qg(t,q), \\
        q(0) &= q_0, \\
        p(T) &= p_1(q(T)). 
    \end{align*}
    The first and third equation are simply an initial value problem in only the $q$ variable (intrinsically, the system on $T^*M$ covers an ODE on $M$ \cite{TrLe2024adj})  and thus, uniquely determines a solution curve $q: [0,T] \rightarrow M$; subsequently, by substituting in this determined curve $q(t)$, the second and fourth equation becomes an ODE in only the $p$ variable with $p(T) = p_1(q(T))$ now fixed and thus, we obtain a unique solution $(q,p): [0,T] \rightarrow T^*M$ (assuming $f$ and $g$ are sufficiently smooth). We will discuss such Hamiltonians in \Cref{sec:adjoint-sensitivity} in the context of adjoint sensitivity analysis.

    We remark that a natural way to define a section of the cotangent bundle $p_1: M \rightarrow T^*M$ is as the exterior derivative of a function $C: M \rightarrow \mathbb{R}$, $p_1(q) = dC(q)$. Interestingly, although all of the Type II d'Alembert variational principles discussed thus far are not generally stationarity conditions, i.e., are not of the form $\delta \mathcal{J} = 0$ for some functional $\mathcal{J}$, the free boundary Type II d'Alembert variational principle is a stationarity condition for such a choice of $p_1 = dC$. To see this, consider the functional on the space of curves $C^1([0,T], T^*M)$ given by 
    $$ \mathcal{C}[q,p] := C(q(T)). $$
        Then, the stationarity of the functional $\mathcal{J}[q,p] := \mathcal{C}[q,p] - S[q,p]$ under partial variations is equivalent to the free boundary Type II d'Alembert variational principle with the choice of $p_1 = dC$, which follows from 
    \begin{align*}
        0 &= \delta \mathcal{J}[q,p] = \delta (\mathcal{C}[q,p] - S[q,p]) = \frac{d}{d}\Big|_{\epsilon = 0} \left(\mathcal{C}[q_\epsilon, p_\epsilon] - \delta S[q_\epsilon,p_\epsilon]\right) \\
        &= \frac{d}{d}\Big|_{\epsilon = 0} C(q_\epsilon(T)) - \delta S[q,p] = \langle dC(q(T)), \delta q(T)\rangle - \delta S[q,p],
    \end{align*}
    which is precisely \eqref{eq:free-type-ii-d'Alembert-vp} with $p_1 = dC$. We will discuss in \Cref{sec:optimization-control} how such a choice of $p_1$ is natural in the contexts of adjoint sensitivity analysis and optimal control.

    Finally, note that the previous Type II d'Alembert variational principles can be thought of as special cases of the free boundary variant. Namely, in the vector space case discussed in \Cref{sec:type-ii-vp-d'Alembert}, the boundary condition $p(T) = p_1$, where $p_1 \in Q^*$ is a fixed cotangent vector, corresponds to a constant section of the bundle $T^*Q \rightarrow Q$ in the free boundary perspective. In this case, since the section is constant, the free boundary problem reduces to a boundary value problem. In the case of a parallelizable manifold discussed in \Cref{sec:type-ii-vp-parallelizable}, the boundary condition $\mu(T) = \mu_1$ corresponds to a constant section of the bundle $M \times V^* \rightarrow M$ and by the inverse of the trivialization, also corresponds to the section of $T^*M$ given by $p_1: q \mapsto \Phi(q)^*\mu_1$. In the trivializing space $M \times V^*$, the section is constant so, as we saw, the associated Hamel system \eqref{eq:model-ham-eom-a}-\eqref{eq:model-ham-eom-d} is a boundary value problem with $\mu(T) = \mu_1$ independent of $q(T)$; however, in canonical coordinates $(q,p)$ on $T^*M$, the system would generally be a free boundary problem with $p(T) = \Phi(q(T))^*\mu_1$ depending on $q(T)$. Thus, for parallelizable manifolds, working in the trivializing space has the advantage of corresponding to a boundary value problem as opposed to a more complicated free boundary problem in canonical coordinates.

\begin{comment}
\textbf{Type II Variational Principles for Maximally Degenerate Hamiltonians.}

\textbf{Type II Variational Principles with an Affine Connection.}
details
\begin{align*}
    S[q,p] = \int_0^T \langle p(t) - \Gamma^*(q)^{T}_t p_f, \dot{q}(t)\rangle - H(t,q(t),p(t)) dt
\end{align*}
IBP
$$ \frac{d}{dt} \Gamma^*(q)^T_t p_f = \nabla^*_{\dot{q}} p_f $$
\end{comment}

\subsubsection{Infinite-dimensional Type II Variational Principles}\label{sec:inf-dim-type-ii-vp}
As a final generalization, we develop here Type II variational principles in the infinite-dimensional setting, namely, on infinite-dimensional reflexive Banach spaces. For a discussion of infinite-dimensional Hamiltonian systems, see \cite{ChMa1974, MaHu1983, MaRa1999}. 

Consider a real reflexive Banach space $Y$ and the associated cotangent bundle with the natural identification $T^*Y \cong Y \times Y^* \ni (\varphi, \pi)$. As in the finite-dimensional case, we denote the duality pairing by $\langle\cdot,\cdot\rangle: Y^* \times Y \rightarrow \mathbb{R}$. Let $\mathfrak{D}^1$ be an affine subspace of $Y$ modeled over a dense subspace $\mathfrak{D}^1_0$ of $Y$. Similarly, let $\mathfrak{D}^2$ be an affine subspace of $Y^*$ modeled over a dense subspace $\mathfrak{D}^2_0$ of $Y^*$. 
\begin{remark}
    Intuitively, we allow the possibility of $\mathfrak{D}^1$ and $\mathfrak{D}^2$ being affine subspaces, rather than the stronger requirement that they are subspaces, in order to allow for inhomogeneous boundary conditions, e.g., when $Y$ corresponds to some function space over a spatial domain. Homogeneous boundary conditions correspond to $\mathfrak{D}^1 = \mathfrak{D}^1_0$ and $\mathfrak{D}^2 = \mathfrak{D}^2_0$. 

    The property that the model subspaces are dense will be used to make the usual fundamental lemma of the calculus of variations argument. Furthermore, this is a typical situation for infinite-dimensional Hamiltonian systems, where these subspaces correspond to the domain of densely defined unbounded linear operators \cite{ChMa1974, MaHu1983, MaRa1999}.
\end{remark}
Then, a Hamiltonian system on the interval $[0,T]$ is specified by a time-dependent mapping which is defined on $\mathfrak{D}^1 \times \mathfrak{D}^2 \subset Y \times Y^*$, i.e.,
$$ \mathcal{H}: [0,T] \times \mathfrak{D}^1 \times \mathfrak{D}^2 \subset [0,T] \times Y \times Y^* \rightarrow \mathbb{R},$$
\sloppy which is furthermore Fr\'{e}chet differentiable with respect to its second and third arguments for all $(t,\varphi,\pi) \in [0,T] \times \mathfrak{D}^1 \times \mathfrak{D}^2$. Note that $D_\varphi \mathcal{H}(t,\varphi,\pi) \in Y^*$ and $D_\pi\mathcal{H}(t,\varphi,\pi) \in Y$. Hamilton's equations associated with this Hamiltonian are
\begin{subequations}\label{eq:inf-dim-Hamilton-eqn}
\begin{align}
    \dot{\varphi}(t) &= D_\pi \mathcal{H}(t,\varphi(t),\pi(t)), \\
    \dot{\pi}(t) &= -D_\varphi \mathcal{H}(t,\varphi(t),\pi(t)).
\end{align}
\end{subequations}
The symplectic structure of such infinite-dimensional Hamiltonian systems is discussed in \cite{ChMa1974}; variational principles for infinite-dimensional Lagrangian systems are discussed in \cite{MaHu1983}. Here, we will instead take a Type II variational perspective. We will first consider an infinite-dimensional generalization of the free boundary Type II variational principle discussed in \Cref{sec:free-boundary-type-ii} and then consider an infinite-dimensional generalization of the Type II variational principle of \citet{LeZh2009}.

Consider the space of curves
$$ \mathfrak{C} := \Big( C^1([0,T],Y) \cap C^0([0,T], \mathfrak{D}^1) \Big) \times \Big(C^1([0,T],Y^*) \cap C^0([0,T], \mathfrak{D}^2)\Big), $$
which we choose in order to make classical sense of solutions to Hamilton's equations \eqref{eq:inf-dim-Hamilton-eqn} (see, e.g., \cite{Ne1992, TrSoLe2024}), although one can in principle use a space of curves in a weaker regularity class for weaker notions of solutions. Note that the spatial boundary conditions for the system are implicitly imposed in the choice of $\mathfrak{D}^1$ and $\mathfrak{D}^2$, see, e.g., \Cref{ex:adjoint-diffusion}. For fixed $\varphi_0 \in \mathfrak{D}^1$, we define the subset of $\mathfrak{C}$ fixing $\varphi(0) = \varphi_0$ to be
$$ \mathfrak{C}_{\varphi_0} := \{ (\varphi, \pi) \in \mathfrak{C} : \varphi(0) = \varphi_0 \}. $$

We define the action $\mathcal{S}: \mathfrak{C}_{\varphi_0} \rightarrow \mathbb{R}$ associated with $\mathcal{H}$ by
$$ \mathcal{S}[\varphi,\pi] = \int_0^T [\langle \pi(t), \dot{\varphi}(t)\rangle - \mathcal{H}(t,\varphi(t),\pi(t))]\,dt. $$
We define the space of partial variations to be
\begin{align*}
    \mathfrak{C}_0 := \Big\{ (\delta\varphi,\delta \pi) \in \Big(& C^1([0,T],Y) \cap C^0([0,T], \mathfrak{D}^1_0) \Big) \\
    & \qquad \times \Big(C^1([0,T],Y^*) \cap C^0([0,T], \mathfrak{D}^2_0)\Big): \delta \varphi(0) = 0 \Big\}.
\end{align*}
Note that $\mathfrak{C}_{\varphi_0} + \mathfrak{C}_0 = \mathfrak{C}_{\varphi_0}$ due to the affine structure of $\mathfrak{D}^1$ and $\mathfrak{D}^2$. Thus, we can define variation of the action to be
\begin{equation}\label{eq:banach-variation}
    \delta \mathcal{S}[\varphi,\pi] := \frac{d}{d\epsilon}\Big|_{\epsilon = 0} \mathcal{S}[\varphi + \epsilon \delta \varphi, \pi + \epsilon \delta \pi]. 
\end{equation}
We define the free boundary Type II d'Alembert variational principle, for a mapping $\pi_1: \mathfrak{D}^1 \rightarrow \mathfrak{D}^2$, to be $\delta \mathcal{S}[\varphi,\pi] = \langle p_1(\varphi(T)),\delta\varphi(T)\rangle$ for all partial variations. Given the regularity of curves in $\mathfrak{C}$, the same argument as in the finite-dimensional case follows through, additionally using the density of $\mathfrak{D}^1_0$ and $\mathfrak{D}^2_0$ to apply the fundamental lemma of the calculus of variations. This variational principle is then equivalent to Hamilton's equations with Type II free boundary conditions:
\begin{subequations}\label{eq:inf-dim-Hamilton-eqn-type-ii}
\begin{alignat}{2}
    \dot{\varphi} &= D_\pi \mathcal{H}(t,\varphi,\pi),&\,&\text{ for all } t\in(0,T),\\
    \dot{\pi} &= -D_\varphi \mathcal{H}(t,\varphi,\pi),&\,&\text{ for all } t\in(0,T), \\
    \varphi(0) &= \varphi_0 \in \mathfrak{D}^1, \\
    \pi(T) &= \pi_1(\varphi(T)) \in \mathfrak{D}^2.
\end{alignat}
\end{subequations}

\begin{example}[Adjoint system for the diffusion equation]\label{ex:adjoint-diffusion}
In \cite{TrSoLe2024}, we explore the Hamiltonian structure associated with adjoint systems for abstract time-dependent semilinear evolution equations on Banach spaces. Here, we will just consider a particular example to contextualize the discussion above. Consider the diffusion equation
$$ \dot{\varphi} = \Delta \varphi, $$
where we take $Y = L^2(\Omega) \cong Y^*$, $\mathfrak{D}^1_0 = \mathfrak{D}^2_0 = \mathfrak{D}^2 = H^2(\Omega)\cap H^1_0(\Omega)$ and $\mathfrak{D}^1 = H^2(\Omega) \cap H^1_g(\Omega)$ where $H^1_g(\Omega)$ is the space of functions in $H^1$ with trace equal to $g \in L^2(\partial \Omega)$, assuming a domain $\Omega \subset \mathbb{R}^n$ whose boundary has sufficient regularity (e.g., Lipschitz) such that the trace $\textup{Tr}: H^1(\Omega) \rightarrow L^2(\partial \Omega)$ is continuous. Define the corresponding adjoint Hamiltonian (in this case, time-independent) to be
$$ \mathcal{H}(\varphi,\pi) = \langle \pi, \Delta \varphi \rangle = \int_\Omega \pi \Delta \varphi\, d^nx. $$
With this Hamiltonian, the system \eqref{eq:inf-dim-Hamilton-eqn-type-ii} arising from the Type II variational principle is the adjoint system for the diffusion equation with the following spatial boundary conditions and Type II (temporal) boundary conditions,
\begin{alignat*}{2}
    \dot{\varphi}(t) &= \Delta \varphi(t),&\,&\text{ for all } t\in(0,T), \\
    \dot{\pi}(t) &= -\Delta \pi(t),&\,&\text{ for all } t\in(0,T), \\
    \textup{Tr}( \varphi(t)) &= g \in L^2(\partial \Omega),&\,&\text{ for all } t\in(0,T), \\
    \textup{Tr}( \pi(t)) &= 0 \in L^2(\partial \Omega),&\,&\text{ for all } t\in(0,T), \\
    \varphi(0) &= \varphi_0 \in \mathfrak{D}^1, \\
    \pi(T) &= \pi_1(\varphi(T)) \in \mathfrak{D}^2.
\end{alignat*}
Note that the initial-boundary value problem for $\varphi$ can be solved in forward time and subsequently, the terminal-boundary value problem for $\pi$ can be solved in reverse time \cite{TrSoLe2024}, since the PDE for $\pi$ is not the diffusion equation but rather the time-reversed diffusion equation (i.e., $-\Delta$ appears on the right hand side of the equation for $\dot{\pi}$ above). The same can not be said for Type III boundary conditions specifying a terminal value for $\varphi$ and an initial value for $\pi$, since the diffusion equation is not well-posed in reverse time. Thus, unlike the finite-dimensional case where well-posedness is preserved under time reversal, this is no longer the case and Type II/III boundary conditions differ in their completeness (in the sense discussed in \Cref{sec:boundary-conditions}), depending on the given problem. Furthermore, Type 0 boundary conditions (i.e., the initial value problem) is not well-posed in this example, since $\pi$ evolves by the time-reversed diffusion equation.
\end{example}

For completeness, we will also extend the Type II variational principle on vector spaces of \citet{LeZh2009} to the infinite-dimensional setting; this is possible since here we are working on Banach spaces instead of more general Banach manifolds. Let $\mathcal{H}$ and its domain be as before. Now, consider the space of curves satisfying Type II boundary conditions $\varphi(0) = \varphi_0 \in \mathfrak{D}^1$ and $\pi(T) = \pi_1 \in \mathfrak{D}^2$, i.e., 
$$ \mathfrak{C}_{\varphi_0, \pi_1} = \{ (\varphi, \pi) \in \mathfrak{C} : \varphi(0) = \varphi_0, \pi(T) = \pi_1 \}, $$
where $\mathfrak{C}$ is as before. We define the space of Type II variations as 
\begin{align*}
\mathfrak{C}_{\textup{II}} &:= \Big\{ (\delta\varphi,\delta \pi) \in \Big( C^1([0,T],Y) \cap C^0([0,T], \mathfrak{D}^1_0) \Big) \\
& \qquad \qquad \times \Big(C^1([0,T],Y^*) \cap C^0([0,T], \mathfrak{D}^2_0)\Big): \delta \varphi(0) = 0, \delta \pi(T) = 0 \Big\}.
\end{align*}
We define the Type II action $\mathcal{S}_{\textup{II}}: \mathfrak{C}_{\varphi_0,\pi_1} \rightarrow \mathbb{R}$ associated with $\mathcal{H}$ by
$$ \mathcal{S}_{\textup{II}}[\varphi,\pi] = \langle \pi(T), \varphi(T)\rangle - \int_0^T [\langle \pi(t), \dot{\varphi}(t)\rangle - \mathcal{H}(t,\varphi(t),\pi(t))]\,dt. $$
Note that $\mathfrak{C}_{\varphi_0, \pi_1} + \mathfrak{C}_{\textup{II}}= \mathfrak{C}_{\varphi_0, \pi_1}$ again due to the affine structure of $\mathfrak{D}^1$ and $\mathfrak{D}^2$. Thus, we can define the variation of the action as before, \Cref{eq:banach-variation}, with Type II variations instead of partial variations. 

We define the Type II variational principle on Banach spaces, for $\varphi_0 \in \mathfrak{D}^1$ and $\pi_1 \in \mathfrak{D}^2$, to be the stationarity of the Type II action $\delta \mathcal{S}_{\textup{II}}[\varphi,\pi] = 0$, for all Type II variations. This variational principle is then equivalent to Hamilton's equations with Type II boundary conditions:
\begin{subequations}\label{eq:evolution-type-ii-pde}
\begin{alignat}{2}
    \dot{\varphi}(t) &= D_\pi \mathcal{H}(t,\varphi,\pi),&\,&\text{ for all } t\in(0,T), \\
    \dot{\pi}(t) &= -D_\varphi \mathcal{H}(t,\varphi,\pi),&\,&\text{ for all } t\in(0,T), \\
    \varphi(0) &= \varphi_0 \in \mathfrak{D}^1, \\
    \pi(T) &= \pi_1 \in \mathfrak{D}^2,
\end{alignat}
\end{subequations}
which again follows from a formal calculation, additionally using the density of $\mathfrak{D}^1_0$ and $\mathfrak{D}^2_0$.

It would be interesting to develop Type II variational integrators in the infinite-dimensional setting, arising from discretizing these Type II variational principles in space and time, and examine their applications, e.g., to adjoint sensitivity analysis of numerical PDEs.

\section{Overview of Hamiltonian Variational Principles}\label{sec:generalizations}
In the previous section, \Cref{sec:hamiltonian-type-ii}, we considered variational principles for Hamiltonian systems with the configuration space being a finite-dimensional vector space $Q$, a finite-dimensional manifold $M$, or a reflexive Banach space $Y$. Particularly, we discussed in detail Type II variational principles for such systems. Now, we provide an overview of generalizations and applications of Hamiltonian variational principles, as well as some connections to the Type II d'Alembert variational principle introduced in the previous section. Namely, we provide an overview of variational principles and applications to variational integrators, optimization and control, Hamiltonian systems with constraints, Hamiltonian PDEs, and stochastic Hamiltonian systems.

\subsection{Lagrangian and Hamiltonian Variational Integrators}\label{sec:variational-integrators}
Geometric numerical integration aims to preserve geometric conservation laws under discretization, and this field is surveyed in the monograph by \citet{HaLuWa2006}. Discrete variational mechanics \cite{MaWe2001, LeOh2008} provides a systematic method of constructing symplectic integrators. It is typically approached from a Lagrangian perspective by introducing the \textit{discrete Lagrangian}, $L_d:Q\times Q\rightarrow \mathbb{R}$, which is a Type I generating function of a symplectic map and approximates the \textit{exact discrete Lagrangian}, which is constructed from the Lagrangian $L: TQ \rightarrow \mathbb{R}$ as
\vspace*{-1.75ex}
\begin{equation}
L_d^E(q_0,q_1;h)=\ext_{\substack{q\in C^2([0,h],Q) \\ q(0)=q_0, q(h)=q_1}} \int_0^h L(q(t), \dot q(t)) dt,\label{exact_Ld_variational}
\end{equation}
which is equivalent to Jacobi's solution of the Hamilton--Jacobi equation.
The exact discrete Lagrangian generates the exact time-$h$ flow map of a Lagrangian system, but, in general, it cannot be computed explicitly. Instead, this can be approximated by replacing the integral with a quadrature formula, and replacing the space of $C^2$ curves with a finite-dimensional function space.

Given a finite-dimensional function space $\mathbb{M}^n([0,h])\subset C^2([0,h],Q)$ and a quadrature formula $\mathcal{G}:C^2([0,h],Q)\rightarrow\mathbb{R}$, $\mathcal{G}(f)=h\sum_{j=1}^m b_j f(c_j h)\approx \int_0^h f(t) dt$, the \textit{Galerkin discrete Lagrangian} is
\begin{align*}
    L_d(q_0,q_1)&=\ext_{\substack{q\in \mathbb{M}^n([0,h]) \\ q(0)=q_0, q(h)=q_1}} \mathcal{G}(L(q, \dot q)) \\
    &=\ext_{\substack{q\in \mathbb{M}^n([0,h]) \\ q(0)=q_0, q(h)=q_1}} h \sum\nolimits_{j=1}^m b_j L(q(c_j h), \dot q(c_j h)).
\end{align*}

Given a discrete Lagrangian $L_d$, the \textit{discrete Hamilton--Pontryagin principle} imposes the discrete second-order condition $q_k^1=q_{k+1}^0$ using Lagrange multipliers $p_{k+1}$, which yields a variational principle on $(Q\times Q)\times_Q T^*Q$,
\begin{equation*}
 \delta \left[\sum\nolimits_{k=0}^{n-1} L_d(q_k^0,q_k^1)+\sum\nolimits_{k=0}^{n-2}p_{k+1}(q_{k+1}^0-q_k^1)\right]=0.
 \end{equation*}
This in turn yields the \textit{implicit discrete Euler--Lagrange equations},
\begin{equation}\label{IDEL} q_k^1=q_{k+1}^0,\qquad p_{k+1}=D_2 L_d(q_k^0, q_k^1), \qquad p_{k}=-D_1 L_d(q_k^0, q_k^1),\end{equation}
where $D_i$ denotes the partial derivative with respect to the $i$-th argument. Making the identification $q_k=q_k^0=q_{k-1}^1$, we obtain the \textit{discrete Lagrangian map} and \textit{discrete Hamiltonian map} which are $F_{L_d}:(q_{k-1},q_k)\mapsto(q_k,q_{k+1})$ and $\tilde{F}_{L_d}:(q_k,p_k)\mapsto(q_{k+1},p_{k+1})$, respectively. The last two equations of \eqref{IDEL} define the \textit{discrete fiber derivatives}, $\mathbb{F}L_d^\pm:Q\times Q\rightarrow T^*Q$,
%\begin{align}
%\newsubeqblock
%\mysubeq \mathbb{F}L_d^+(q_k, q_{k+1})&=(q_{k+1},D_2 L_d(q_k, q_{k+1})), \label{FLd+}\\
%\mysubeq \mathbb{F}L_d^-(q_k, q_{k+1})&=(q_k,-D_1 L_d(q_k, q_{k+1})).\label{FLd-}
%\end{align}
\begin{align*}
\mathbb{F}L_d^+(q_k, q_{k+1})&=(q_{k+1},D_2 L_d(q_k, q_{k+1})),\\
\mathbb{F}L_d^-(q_k, q_{k+1})&=(q_k,-D_1 L_d(q_k, q_{k+1})).
\end{align*}
These two discrete fiber derivatives induce a single unique \textit{discrete symplectic form} $\Omega_{L_d}=(\mathbb{F}L_d^\pm)^*\Omega$, where $\Omega$ is the canonical symplectic form on $T^*Q$, and the discrete Lagrangian and Hamiltonian maps preserve $\Omega_{L_d}$ and $\Omega$, respectively.
%%%%%%%%%%%%%%
The discrete Lagrangian and Hamiltonian maps can be expressed as $F_{L_d}=(\mathbb{F}L_d^-)^{-1}\circ \mathbb{F}L_d^+$ and $\tilde{F}_{L_d}=\mathbb{F}L_d^+\circ (\mathbb{F}L_d^-)^{-1}$, respectively. This characterization allows one to relate the approximation error of the discrete flow maps to the approximation error of the discrete Lagrangian.

The variational integrator approach simplifies the numerical analysis of symplectic integrators. The task of establishing the geometric conservation properties and order of accuracy of the discrete Lagrangian map $F_{L_d}$ and discrete Hamiltonian map $\tilde{F}_{L_d}$ reduces to the simpler task of verifying certain properties of the discrete Lagrangian $L_d$ instead.

\begin{theorem}[Discrete Noether's theorem (Theorem 1.3.3 of \cite{MaWe2001})] If a discrete Lagrangian $L_d$ is invariant under the diagonal action of $G$ on $Q\times Q$, then the single unique \textit{discrete momentum map}, $\mathbf{J}_{L_d}=(\mathbb{F}L_d^\pm)^*\mathbf{J}$, is invariant under the discrete Lagrangian map $F_{L_d}$, i.e., $F_{L_d}^* \mathbf{J}_{L_d}=\mathbf{J}_{L_d}$.
\end{theorem}

\begin{theorem}[Variational error analysis (Theorem 2.3.1 of \cite{MaWe2001}; see also \cite{PaCu2009})\label{thm_variational_error_analysis}]
If a discrete Lagrangian $L_d$ approximates the exact discrete Lagrangian $L_d^E$ to order $p$, i.e., $ L_d(q_0, q_1;h)=L_d^E(q_0,q_1;h)+\mathcal{O}(h^{p+1}),$
then the discrete Hamiltonian map $\tilde{F}_{L_d}$ is an order $p$ accurate one-step method.
\end{theorem}

The bounded energy error of variational integrators can be understood by performing backward error analysis, which then shows that the discrete flow map is approximated with exponential accuracy by the exact flow map of the Hamiltonian vector field of a modified Hamiltonian~\cite{BeGi1994,Ta1994}. For examples, applications, and generalizations of such Lagrangian variational integrators, see \cite{MaWe2001, Ob2016, deLe2007, FaLeLeMcSc2006, LeLeMc2005a, LeLeMc2007b, LeMaOr2008, LeMaOrWe2004, HaLu2020, XuTrLe2022, HaLe2015, ObSa2015, ZaNeRoSc2019, Stern2015, ZhQiBuBh2014}.

Given a degenerate Hamiltonian, where the Legendre transform $\mathbb{F}H:T^*Q\rightarrow TQ$, $(q,p)\mapsto (q, \frac{\partial H}{\partial p})$, is noninvertible, there is no equivalent Lagrangian formulation. Thus, a characterization of variational integrators directly in terms of the continuous Hamiltonian is desirable. This is achieved by considering the Type II analogue of Jacobi's solution \cite{LeZh2009}, given by
\[H_d^{+,E}(q_k, p_{k+1})=  \ext_{\substack{(q, p) \in
C^1([t_k,t_{k+1}],T^*Q)\\q(t_k)=q_k, p(t_{k+1})=p_{k+1}}} \Big[ p(t_{k+1}) q (t_{k+1}) - \int_{t_k}^{t_{k+1}} \left[ p \dot{q}-H(q, p) \right]
dt \Big].
\]
A computable Galerkin discrete Hamiltonian $H_d^+$ is obtained by choosing a finite-dimensional function space and a quadrature formula,
\begin{align*}
    H_d^+&(q_0,p_1)\\
    &=\ext_{\substack{q\in \mathbb{M}^n([0,h]) \\ q(0)=q_0\\(q(c_j h), p(c_j h))\in T^*Q}} \left[p_1 q(t_1) - h \sum\nolimits_{j=1}^m b_j [ p(c_j h)\dot q(c_j h)-H(q(c_j h),p(c_j h)) ]\right].
\end{align*}
Interestingly, the Galerkin discrete Hamiltonian does not require a choice of a finite-dimensional function space for the curves in the momentum, as the quadrature approximation of the action integral only depend on the momentum values $p(c_j h)$ at the quadrature points, which are determined by the extremization principle. In essence, this is because the action integral does not depend on the time derivative of the momentum $\dot{p}$. As such, both the Galerkin discrete Lagrangian and the Galerkin discrete Hamiltonian depend only on the choice of a finite-dimensional function space for curves in the position, and a quadrature rule. It was shown in Proposition 4.1 of \citep{LeZh2009} that when the Hamiltonian is hyperregular, and for the same choice of function space and quadrature rule, they induce equivalent Galerkin variational integrators. However, as shown in \cite{ScLe2017}, for more general approximation schemes, discrete Lagrangian and Hamiltionian variational integrators need not agree, even assuming hyperregularity; this provides another reason to study Hamiltonian variational integrators directly from Hamiltonian variational principles.

The \textit{Type II discrete Hamilton's phase space variational principle} states that 
\[
\delta \left\{p_N q_N -\sum_{k=0}^{N-1}\left[p_{k+1} q_{k+1} -H_d^+(q_k,
 p_{k+1})\right]\right\}=0,
\]
for discrete curves in $T^*Q$ with fixed $(q_0, p_N)$ boundary conditions. This yields the \textit{discrete Hamilton's equations}, which are given by
\begin{equation}\label{dHamiltonEq}
q_{k+1} = D_2 H_d^+(q_k, p_{k+1}), \qquad p_k=D_1 H_d^+(q_k, p_{k+1}).
\end{equation}

Given a discrete Hamiltonian $H_d^+$, we introduce the \textit{discrete fiber derivatives} (or discrete Legendre transforms), $\mathbb{F}^+ H_d^+$,
\begin{align*}
\mathbb{F}^+H_d^+&: (q_0,p_1)\mapsto (D_2 H_d^+(q_0,p_1),p_1),\\
\mathbb{F}^-H_d^+&: (q_0,p_1)\mapsto (q_0 , D_1 H_d^+ (q_0,p_1)).
\end{align*}
The discrete Hamiltonian map can be expressed in terms of the discrete fiber derivatives,
$$ \tilde{F}_{H_d^+}(q_0,p_0) = \mathbb{F}^+H_d^{+} \circ (\mathbb{F}^-H_d^+)^{-1} (q_0,p_0) = (q_1,p_1) ,$$

Similar to the Lagrangian case, we have a discrete Noether's theorem and variational error analysis result for Hamiltonian variational integrators.

\begin{theorem}[Discrete Noether's theorem (Theorem 5.3 of \cite{LeZh2009})] Given the action $\Phi$ on the configuration manifold $Q$, let $\Phi^{T^*Q}$ be the cotangent lifted action on $T^*Q$.  If the generalized discrete Lagrangian $R_d(q_0, q_1, p_1)=p_1 q_1 -H_d^+(q_0, p_1)$ is invariant under the cotangent lifted action $\Phi^{T^*Q}$, then the discrete Hamiltonian map $\tilde{F}_{H_d^+}$ preserves the momentum map, i.e., $\tilde{F}_{H_d^+}^* \mathbf{J}=\mathbf{J}$.
\end{theorem}

\begin{theorem}[Variational error analysis (Theorem 2.2 of \cite{ScLe2017})]
If a discrete Hamiltonian $H_d^+$ approximates the exact discrete Hamiltonian $H_d^{+,E}$ to order $p$, i.e., $ H_d^+(q_0, p_1;h)=H_d^{+,E}(q_0,p_1;h)+\mathcal{O}(h^{p+1})$, then the discrete Hamiltonian map $\tilde{F}_{H_d^+}$ is an order $p$ accurate one-step method.
\end{theorem}

It should be noted that there is an analogous theory of discrete Hamiltonian variational integrators based on Type III generating functions $H_d^-(p_0,q_1)$. 

Hamiltonian variational integrators were further generalized to the stochastic setting \cite{HoTy2018, KrTy2020} (see \Cref{sec:stochastic}), the presymplectic setting \cite{TrLe2024adj} and the setting of multisymplectic PDEs \cite{TrLe2022} (see \Cref{sec:Hamiltonian-PDEs}), based upon discretizations of the respective Type II variational principles. Lie--Poisson Hamiltonian variational integrators are discussed in \cite{MaRo2010, deDiego2018}. Hamiltonian variational integrators for optimal control are discussed in \cite{deLe2007}. The application of Hamiltonian variational integrators to symplectic accelerated optimization is examined in \cite{SAO1,SAO2,SAO4,SAO5,SAO6}.

\subsection{Optimization and Control}\label{sec:optimization-control}
In the following sections, we discuss Hamiltonian variational principles in the context of optimization and control; particularly, we discuss applications in adjoint sensitivity analysis, optimal control and symplectic accelerated optimization. 

\subsubsection{Adjoint Sensitivity Analysis}\label{sec:adjoint-sensitivity}
The solution of many nonlinear problems involves successive linearization, and as such variational equations and their adjoints play a critical role in a variety of applications. Adjoint equations are of particular interest when the parameter space is significantly higher dimension than that of the output or objective. In particular, the simulation of adjoint equations arise in sensitivity analysis~\cite{Ca1981, CaLiPeSe2003}, adaptive mesh refinement~\cite{LiPe2003}, uncertainty quantification~\cite{WaDuAlIa2012}, automatic differentiation~\cite{Gr2003}, superconvergent functional recovery~\cite{PiGi2000}, optimal control~\cite{Ro2005}, optimal design~\cite{GiPi2000}, optimal estimation~\cite{NgGeBe2016}, and deep learning viewed as an optimal control problem~\cite{DeCeEhOwSc2019}.

The study of geometric aspects of adjoint systems arose from the observation that the combination of any system of differential equations and its adjoint equations are described by a formal Lagrangian~\cite{Ib2006, Ib2007}. This naturally leads to the question of when the formation of adjoints and discretization commutes~\cite{SoTz1997}, and prior work on this include the Ross--Fahroo lemma~\cite{RoFa2001}, and the observation by \citet{Sa2016} that the adjoints and discretization commute if and only if the discretization is symplectic, in the context of Runge--Kutta time discretization. In \cite{TrLe2024adj}, we studied adjoint systems associated with ODEs and differential-algebraic equations (DAEs), and their discretization, using techniques from symplectic geometry, presymplectic geometry and variational principles. In \cite{TrSoLe2024}, we extended this study to the infinite-dimensional setting of semilinear evolution equations on Banach spaces. Understanding the geometry of adjoint systems led us to a general classification of the commutativity of discretize-then-optimize versus optimize-then-discretize approaches to adjoint sensitivity analysis \cite{TrLe2024adj, TrSoLe2024}.

Here, we discuss the free boundary Type II d'Alembert variational principle (\Cref{sec:free-boundary-type-ii}) in the context of adjoint sensitivity analysis, where we show that it has a natural interpretation in this setting.

Let the cost function be
\begin{equation}\label{eq:terminal-and-running-cost}
    C(q(T)) + \int_0^T g(t,q(t))\, dt,
\end{equation}
where $C$ is a terminal cost function and the above integral of $g$ is a running cost, subject to the dynamics of an ODE $\dot{q}(t) = f(t,q(t)), q(0)=q_0$ on a manifold $M$. Consider the problem of finding the sensitivity of this cost function with respect to variations in the initial condition, i.e., we wish to compute the derivative of the cost with respect to the initial condition $q(0)$, viewing the cost function implicitly as a function of the initial condition. The standard adjoint approach for computing this sensitivity \cite{Ca1981, CaLiPeSe2003} is to consider the augmented cost functional
$$ \mathcal{J}[q,p] := C(q(T)) + \int_0^T g(t,q(t))\, dt - \int_0^T \langle p(t), \dot{q}(t) - f(t,q(t))\rangle\,dt. $$
Here, the curve $p(t)$ is interpreted as a Lagrange multiplier which enforces the ODE constraint $\dot{q}=f$. Taking the variation of $\mathcal{J}$ with respect to arbitrary $\delta q$ and $\delta p$ shows that if one chooses $p(T) = dC(q(T))$ and $p(t)$ to satisfy the augmented adjoint equation
$$ \dot{p} = - [D_qf(t,q)]^*p - D_q g(t,q), $$
then the desired sensitivity is given by $p(0)$.

We will reinterpret this in the context of the free boundary Type II d'Alembert principle. We define the augmented adjoint Hamiltonian associated with this cost function as
\begin{align}
    H_g&: [0,T] \times T^*M \rightarrow \mathbb{R}, \label{eq:adjoint-augmented-Hamiltonian} \\
    H_g(t,q,p) &= \langle p,f(t,q)\rangle + g(t,q). \nonumber
\end{align}
The associated augmented action is
\begin{equation}\label{eq:adjoint-augmented-action}
    S_g[q,p] = \int_0^T [\langle p(t), \dot{q}(t)\rangle - H_g(t,q(t),p(t))]\, dt.
\end{equation}
Now, we apply the free boundary Type II d'Alembert variational principle \eqref{eq:free-type-ii-d'Alembert-vp} with $p_1 = dC$, i.e.,
$$ \delta S_g[q,p] = \langle dC(q(T)), \delta q(T)\rangle, $$
where $\delta q$ satisfies $\delta q(0)=0$ and $\delta p$ is arbitrary. This yields the augmented adjoint system
\begin{subequations}\label{eq:augmented-adjoint-system-type-ii}
\begin{align}
    \dot{q}(t) &= f(t,q(t)),\quad q(0) = q_0,\\
    \dot{p}(t) &= - [D_qf(t,q(t))]^*p(t) - D_q g(t,q(t)), \quad p(T) = dC(q(T)).
\end{align}
\end{subequations}
As we noted in \Cref{sec:free-boundary-type-ii}, this variational principle can be expressed as the stationarity condition
$$ \delta (\mathcal{C}[q,p] - S_g[q,p]) = 0, $$ 
for all partial variations, where $\mathcal{C}[q,p]=C(q(T))$. Noting that $\mathcal{C}[q,p] - S_g[q,p]$ is precisely the augmented cost functional $\mathcal{J}[q,p]$, the free boundary Type II d'Alembert variational principle is precisely $\delta \mathcal{J}=0$
subject to partial variations. Thus, in the setting of adjoint sensitivity analysis, the free boundary Type II d'Alembert variational principle is naturally interpreted as the stationarity of the augmented cost functional. Finally, to obtain the desired sensitivity, we compute the variation of $\mathcal{J}$ with $(q,p)$ satisfying \eqref{eq:augmented-adjoint-system-type-ii} without imposing $\delta q(0) = 0$. For such variations, $\mathcal{J}$ is no longer stationary; instead, by similar calculations to \Cref{sec:type-ii-vp-d'Alembert}, we have
$$ \delta \mathcal{J}[q,p] = \langle p(0), \delta q(0)\rangle, $$
which again shows that the desired sensitivity is $p(0)$.

\subsubsection{Optimal Control}\label{sec:optimal-control}
The Type II d'Alembert variational principle can also be extended to the setting of optimal control; as we will see, this variational principle will result in the optimality conditions in Pontryagin's maximum principle \cite{BuLe2014, Bloch2015, deLe2007}. Consider the following fixed-time optimal control problem with state dynamics on $M$,
\begin{align}\label{eq:OCP}
        \min_{u:[0,T]\rightarrow U} &\ C(q(T)) + \int_0^T g(t,q(t),u(t))\,dt \\
        &\textbf{subject to }  \nonumber \\
            & \quad \dot{q}(t) = f(t,q(t),u(t)), t \in (0,T),\nonumber \\
            & \quad q(0)=q_0. \nonumber
\end{align}
That is, the problem is to find a control $u(t)$, i.e., a curve on the admissible control set $U$ (which we assume to be a smooth manifold), such that the cost function is minimized, subject to the initial value state dynamics on $M$, where the running cost function $g$ and the state dynamics $f$ additionally depend on the control $u$.

Let us apply the free boundary Type II d'Alembert variational principle with the section $p_1 = dC$. We define the augmented Hamiltonian and augmented action in analogy with \eqref{eq:adjoint-augmented-Hamiltonian} and \eqref{eq:adjoint-augmented-action} where we now include the dependence on $u$,
\begin{align*}
    H_g &: [0,T] \times T^*M \times U \rightarrow \mathbb{R}, \\
    H_g(t,q,p,u) &= \langle p, f(t,q,u)\rangle + g(t,q,u), \\ 
    S_g[q,p,u] &= \int_0^T [\langle p(t), \dot{q}(t)\rangle - H_g(t,q(t),p(t),u(t))]\, dt.
\end{align*}
\sloppy A direct calculation shows that the free boundary Type II d'Alembert variational principle $\delta S_g[q,p,u] = \langle dC(q(T)), \delta q(T)\rangle $ for all partial variations, i.e., variations $(\delta q,\delta p, \delta u)$ such that $\delta q(0)=0$, yields the following system:
\begin{subequations}
    \begin{alignat}{2}
        \dot{q}(t) &= f(t,q(t), u^*(t)),&\,&\text { for all }t\in(0,T), \label{eq:max-principle-a} \\
        \dot{p}(t) &= - [D_qf(t,q(t),u^*(t))]^*p(t) - D_q g(t,q(t)),&\,&\text { for all }t\in(0,T), \label{eq:max-principle-b} \\
        u^*(t) &= \text{arg\ ext}_{u \in U} H(t,q(t),p(t), u),&\,&\text { for all }t\in(0,T), \label{eq:max-principle-c} \\
        q(0) &= q_0,\label{eq:max-principle-d} \\
        p(T) &= dC(q(T)). \label{eq:max-principle-e}
    \end{alignat}
\end{subequations}
Equations \eqref{eq:max-principle-a}-\eqref{eq:max-principle-e} are the well-known necessary conditions for optimality of Pontryagin's maximum principle \cite{BuLe2014, Bloch2015, deLe2007}. Equations \eqref{eq:max-principle-a}-\eqref{eq:max-principle-b} are the Hamiltonian dynamics of an optimal trajectory; equation \eqref{eq:max-principle-c} corresponds to stationarity with respect to variations $\delta u$; equation \eqref{eq:max-principle-d} and \eqref{eq:max-principle-e} are the transversality conditions for an optimal trajectory corresponding to essential and natural boundary conditions from the variational principle, respectively.

For a discussion of the use of variational integrators for solving the free boundary system \eqref{eq:max-principle-a}-\eqref{eq:max-principle-e}, see \cite{ChHaVi2009, CoGuBlMa2015, JiKoMa2013, Ob2016, deLe2007}. For further discussion of variational principles in the context of optimal control, see \cite{MaZo2021, Pe2022, deLe2007, Bloch2015, Hest1965, Tro1996, vp4}.

\subsubsection{Symplectic Accelerated Optimization}
Nesterov's accelerated gradient method was introduced in \cite{Nes83}, and was shown to converge in $\mathcal{O}(1/k^2)$ to the minimum of the convex objective function $f$, improving on the $\mathcal{O}(1/k)$ convergence rate exhibited by the standard gradient descent methods.
This $\mathcal{O}(1/k^2)$ convergence rate was shown in \cite{Nes04} to be optimal among first-order methods using only information about $\nabla f$ at consecutive iterates. This phenomenon in which an algorithm displays this improved rate of convergence is referred to as acceleration, and it was shown in \cite{SuBoCa16} that Nesterov's accelerated gradient method limits to a second-order ODE, as the time-step goes to 0, and that the objective function $f(x(t))$ converges to its optimal value at a rate of $\mathcal{O}(1/t^2)$ along the trajectories of this ODE.

Symplectic optimization generalizes Nesterov's accelerated gradient method by systematically discretizing the flow of the $p$-parameterized family of Bregman Hamiltonians~\cite{WiWiJo2016},
\begin{equation} 
	H(x,r,t)= pt^{p-1} \left[   D_{h^*}(\nabla h(x)+t^p r , \nabla h(x) ) + Ct^pf(x)  \right],
\end{equation}
whose flow converges to the minimizer of $f(x)$ at the rate $\mathcal{O}(1/t^p)$. Since the Hamiltonian is time-dependent, direct discretization using symplectic integrators is potentially problematic, as such methods are generally developed for time-independent problems. This is addressed using a Poincar\'e transformed Hamiltonian~\cite{ZaSz1975}. More precisely, given a Hamiltonian $H(q,p)$, and a desired transformation of time $t \mapsto \tau$ described by the monitor function $g(q,p)$ via
\begin{equation}
	\frac{dt}{d\tau} = g(q,p),
\end{equation}
there is a Poincar\'e transformed Hamiltonian system on extended phase space,
\begin{equation} \label{TransformedH}
	\bar{H}(\bar{q},\bar{p}) = g(q,p) \left(H(q,p) + p^t \right),
\end{equation}
where $\bar{q} = \begin{bmatrix} q \\ q^t \end{bmatrix} $ and $\bar{p} = \begin{bmatrix} p \\ p^t \end{bmatrix} $. One commonly chooses $q^t=t$ and $p^t=-H(q(0),p(0))$, so that $\bar{H}(\bar{q},\bar{p})=0$ along all integral curves through $(\bar{q}(0),\bar{p}(0))$. The time $t$ is the physical time, while $\tau$ is the fictive time. The Poincar\'e transformation was combined with Hamiltonian variational integrators~\cite{LeZh2009} to obtain time-adaptive Hamiltonian variational integrators~\cite{SAO1}, which can be applied to symplectic optimization as follows. Since the $p$-Bregman Hamiltonians exhibit a time-dilation symmetry with time-rescaling $\frac{dt}{d\tau} = \frac{p}{\mathring{p}} t^{1-\frac{\mathring{p}}{p}}$, the Poincar\'e transformed Hamiltonian for this rescaling is given by,
\begin{equation} \label{ADHamiltonian}
	\bar{H}(\bar{q},\bar{r}) =  \frac{1}{\mathring{p}} \left[ \frac{p^2}{2(q^t)^{p+\frac{\mathring{p}}{p}} }  \langle r , r \rangle  + Cp^2 (q^t)^{2p-\frac{\mathring{p}}{p}} f(q) + p r^t (q^t)^{1-\frac{\mathring{p}}{p}}    \right] . 
\end{equation}
The resulting time-adaptive Hamiltonian discretization~\cite{SAO1} yields a significantly more efficient optimization algorithm than direct discretization, and the application to the Bregman Hamiltonians outperform Nesterov's accelerated gradient method.

A generalization of the Bregman Hamiltonian to Riemannian manifolds was introduced in \cite{SAO2}. For geodesically convex and weakly-quasi-convex cost functions, the Hamiltonian is given by
\begin{equation} \label{pBregmanH} H_{p}(X,R,t)= \frac{p}{2t^{\lambda^{-1}\zeta p +1}} \langle\!\langle R , R\rangle\!\rangle + Cpt^{(\lambda^{-1}\zeta +1)p-1} f(X).	
\end{equation}
The corresponding Hamiltonian flow converges to the minimizer at the rate of $\mathcal{O}(1/t^p)$. For the geodesically $\mu$-strongly convex case, the Hamiltonian is given by
\begin{equation} \label{BregmanHGeneralSC}
	H^{SC}(X,R,t)= \frac{e^{- \eta t} }{2} \langle\!\langle R , R\rangle\!\rangle + e^{\eta t} f(X),
\end{equation}
where $\eta = \left(\frac{1}{\sqrt{\zeta}} + \sqrt{\zeta} \right) \sqrt{\mu}$, and $\zeta$ depends on the lower bound of the sectional curvature of the Riemannian manifold and the diameter of the domain. The flow converges to the minimizer at the rate of $\mathcal{O}(e^{-\sqrt{\frac{\mu}{\zeta}}t})$. Discretizing these flows yield efficient Riemannian optimization algorithms on embedded manifolds, but the challenge of developing symplectic optimization algorithms intrinsically on manifolds and their discrete-time convergence analysis remains, and the intrinsic Type II variational principles discussed in this paper may provide insights to this problem. In particular, the analogue of Type II generating functions associated with such variational principles will be relevant to the construction of Hamiltonian variational integrators on manifolds.

\subsection{Hamiltonian Partial Differential Equations}\label{sec:Hamiltonian-PDEs}
We consider Hamiltonian PDEs as a generalization of Hamiltonian ODEs, incorporating additional spatial dependence into the phase space variables. Hamiltonian PDEs can be approached from an evolution equation perspective or from a multisymplectic perspective. Other formulations exist, such as polysymplectic formulations, but we will consider only these two.

In the evolution equation perspective of a PDE, one interprets the state configuration or \textit{field} as varying in space and evolving in time (there are many works on PDE, see, e.g., the classic text \cite{Ev1998}). In this setting, Hamiltonian PDEs are infinite-dimensional Hamiltonian systems evolving on an appropriate Banach space or Banach manifold; see, e.g., \cite{ChMa1974, MaMoRa1985, Struwe2008, MaHu1983, MaRa1999}. We have already discussed an infinite-dimensional Type II variational principle from the evolution equation perspective in \Cref{sec:inf-dim-type-ii-vp}, so we will move on to discuss the multisymplectic perspective.

Instead of viewing a Hamiltonian PDE as an infinite-dimensional Hamiltonian system with time as the independent variable, in the multisymplectic formulation, spacetime coordinates are the independent variables, which has been extensively studied in, for example, \cite{GoIsMaMo1998, GoIsMaMo2004, MaSh1999, MaPeShWe2001, MaPaSh1998}. The multisymplectic formulation is a manifestly covariant formulation of the equations of motion, and admits covariant analogues of symplecticity and Noether's theorem. 

In \citet{VaLiLe2011}, a Type II variational principle for multisymplectic Hamiltonian PDEs is introduced, which we outline here. Consider a trivial vector bundle $E = X \times Q \rightarrow X$ over an oriented spacetime $X$ with boundary $\partial X$, with volume form denoted by $d^{n+1}x$. Let $\Theta$ be the Cartan form on the dual jet bundle $J^1E^*$, which has coordinates $(x^\mu, \phi^A, p, p^A_\mu)$, where $x^\mu$ are the coordinates on spacetime, $\phi^A$ are the coordinates on $Q$, and $p$ and $p^A_\mu$ are the coordinates of the affine map on the jet bundle, $v^A_\mu \mapsto (p + p_A^\mu v^A_\mu)d^{n+1}x$. Define the restricted dual jet bundle $\widetilde{J^1E^*}$ as the quotient of $J^1E^*$ by horizontal one-forms; this space is coordinatized by $(x^\mu,\phi^A,p^A_\mu)$ and is the relevant configuration bundle for a multisymplectic Hamiltonian PDE; we interpret $\phi^A$ as the value of the field and $p_A^\mu$ as the associated momenta in the direction $x^\mu$. The dual jet bundle can be viewed as a bundle over the restricted bundle, $\mu : J^1E^* \rightarrow \widetilde{J^1E^*}$ (see \citet{LePrRoVi2017}). In this setting, the Hamiltonian is a map $H:\widetilde{J^1E^*} \rightarrow \mathbb{R}$. This defines a section of the bundle $\mu$, in coordinates $\tilde{H}(x^\mu,\phi^A,p_A^\mu) = (x^\mu, \phi^A, -H, p_A^\mu)$ or using the projections $\pi^{j,k}$ from the bundle of $(j+k)$-forms on $E$ to the subbundle of $j$-horizontal, $k$-vertical forms, this can be defined as the set of $z \in J^1E^*$ such that $\pi^{n+1,0}(z) = - H(\pi^{n,1}(z))d^{n+1}x$. Using this section, one can pullback the Cartan form to a form on the restricted bundle,
$$ \Theta_H = \tilde{H}^*\Theta = p_A^\mu d\phi^A \wedge d^nx_\mu - H d^{n+1}x. $$
The analogue of Type II boundary conditions for Hamiltonian PDEs is to consider a partition of the boundary $\partial X = (\partial X)_0 \cup (\partial X)_1$ where $(\partial X)_0 \cap (\partial X)_1 = \emptyset$. The Type II boundary conditions are given by specifying the values of the field $\phi$ on $(\partial X)_0$ and the values of the normal momenta, denoted $p^n$, on $(\partial X)_1$, i.e., 
\begin{equation}\label{eq:pde-type-ii-bc}
    \phi|_{(\partial X)_0} = \varphi_0,\quad p^n|_{(\partial X)_1} = \pi_1,
\end{equation}
for given boundary data $\varphi_0$ and $\pi_1$.
Define the Type II action $S$ as a functional on the sections of $\widetilde{J^1E^*} \rightarrow X$ satisfying the Type II boundary conditions \eqref{eq:pde-type-ii-bc},
\begin{equation}\label{eq:PDE-type-ii-action}
S[\phi,p] = \int_{(\partial X)_1} p^\mu \phi d^{n}x_\mu - \int_X (\phi,p)^* \Theta_H.
\end{equation}
The Type II variational principle states that this action is stationary for variations which preserve the boundary conditions \eqref{eq:pde-type-ii-bc}, i.e., $\delta \phi|_{(\partial X)_0} = 0$ and $\delta p^n|_{(\partial X)_1}=0$. This yields
$$ 0 = \delta S[\phi,p] = \int_U(\phi,p)^* i_V d\Theta_H, $$
where we defined the multisymplectic form $\Omega_H = -d\Theta_H$ and $V$ denotes the variation $V := \delta\phi^A \partial/\partial \phi^A + \delta p_A^\mu \partial/\partial p_A^\mu$. The equations of motion associated with this variational principle are the De Donder--Weyl equations with Type II boundary conditions,
\begin{subequations}
    \begin{align}
        \partial_\mu\phi^A &= \frac{\partial H}{\partial p_A^\mu}, \\
        \partial_\mu p_A^\mu &= - \frac{\partial H}{\partial \phi^A},\\
        \phi|_{(\partial X)_0} &= \varphi_0, \\
        p^n|_{(\partial X)_1} &= \pi_1,
    \end{align}
\end{subequations}
which can be thought of as a covariant analogue of \Cref{eq:evolution-type-ii-pde}. Associated with multisymplectic Hamiltonian PDEs is a multisymplectic conservation law which is a covariant generalization of symplecticity. This is given by the multisymplectic form formula
$$ \int_{\partial U} (\phi,p)^*(i_V i_W \Omega_H) = 0, $$
where $V$ and $W$ are any first variations of De Donder--Weyl equations. 

Multisymplectic integrators for Hamiltonian PDEs are numerical methods which preserve at the discrete level this multisymplectic conservation law, see \cite{IsSc2006, BrRe2001, BrRe2006, BrRe2001b, McAr2020, HoLiSu2006, RyMcFr2007, Re2000b, IsSc2004, Re2000, MaPaSh1998}. In \cite{TrLe2022}, we introduce a variational construction of multisymplectic integrators based on a finite element discretization and quadrature approximation of the above Type II variational principle. This provides a systematic method to construct multisymplectic integrators based on a chosen finite element space and quadrature rule; for example, we show that nodal finite element spaces and nodal quadrature rules produce multisymplectic Runge--Kutta methods \cite{HoLiSu2006, RyMcFr2007, Re2000}, although more general approximations are possible within this construction.

\subsection{Constrained Hamiltonian Systems}\label{sec:constrained-hamiltonian-systems}
There are many variational principles associated with constrained Lagrangian and Hamiltonian systems (see, for example, \cite{Bloch2015, vp1, vp2, vp4, vp5, vp6}); such variational principles arise, for example, in the control of nonholonomic systems \cite{Bloch2015, vp4}. Here, we will consider the Hamilton--d'Alembert principle in phase space as developed in \cite{dirac1,dirac2}. Note that discrete analogues of the Dirac structures and associated variational principles discussed in \cite{dirac1,dirac2} were developed in \cite{dirac3}.

For this discussion, we will work in coordinates, although this variational principle was also developed intrinsically in \cite{dirac1, dirac2} via Dirac structures. Let $\Delta_M \subset TM$ be a regular distribution on $M$ (see, e.g., \cite{dgtext}) and denote the corresponding subspace of $T_qM$ at $q \in M$ as $\Delta(q)$. Let $\Delta^\circ_M \subset T^*M$ denote the annihilator of $\Delta_M$ and similarly denote the corresponding subspace of $T^*M$ at $q \in M$ as $\Delta^\circ(q) \subset T^*_qM$. Define the space of curves 
\begin{align*}
    C^1_{\Delta_M, q_0, q_1}&([0,T], T^*M) \\
    &:= \{ (q,p) \in C^1([0,T],T^*M): \dot{q}(\cdot) \in \Delta(q(\cdot)), q(0)=q_0, q(T)=q_1\},
\end{align*} 
and the action
\begin{align*}
    &S: C^1_{\Delta_M, q_0, q_1}([0,T], T^*M)  \rightarrow \mathbb{R}, \\
    &S[q,p] := \int_0^T [\langle p,\dot{q}\rangle - H(q,p)]\, dt.
\end{align*}
The Hamilton--d'Alembert principle in phase space \cite{dirac2} is then given by the stationarity condition 
$$ \delta S[q,p] = 0, $$
for all variations $(\delta q, \delta p)$ such that $\delta q(\cdot) \in \Delta (q(\cdot))$ and $\delta q(0) = 0$, $\delta q(T) = 0$, corresponding to fixed position endpoints. The variation is computed to be
$$ \delta S[q,p] = \langle p(T), \delta q(T)\rangle - \langle p(0), \delta q(0)\rangle  + \int _0^T \left[\langle \delta p(t), \dot{q}(t) - D_p H \rangle + \langle -\dot{p} - D_qH, \delta q \rangle \right]\,dt. $$
This variational principle yields the nonholonomic equations in Hamiltonian form, with fixed endpoint conditions, i.e.,
\begin{alignat*}{2}
    \dot{q}(t) &= D_pH(q(t),p(t)) \in \Delta(q(t)),&\,&\text { for all }t\in(0,T), \\
    \dot{p}(t) &+ D_qH(q(t),p(t)) \in \Delta^\circ(q(t)),&\,&\text { for all }t\in(0,T), \\
    q(0) &= q_0, \\
    q(T) &= q_1.
\end{alignat*}
However, as recently observed in \cite{vp3}, completely fixing the boundary conditions for constrained systems may impose too many conditions and thus, the variational principle may generically give no solutions. In light of this, the conditions $\delta q(0) = 0$ and $\delta q(T) = 0$ may be too strong given that $\delta q(\cdot)$ is not arbitrary but rather $\delta q(\cdot) \in \Delta(q(\cdot))$. Let us, for example, relax the condition at the terminal time $\delta q(T) = 0$ to instead be $\delta q(T) \in \Delta (q(T))$, while still imposing $\delta q(0) = 0$. The variation of the action is now
$$ \delta S[q,p] = \langle p(T), \delta q(T)\rangle + \int _0^T \left[\langle \delta p(t), \dot{q}(t) - D_p H \rangle + \langle -\dot{p} - D_qH, \delta q \rangle \right]\,dt. $$
Demanding stationarity of the action then gives a natural boundary condition $p(T) \in \Delta^\circ(q(T))$, since otherwise, if $p(T) \not\in \Delta^\circ(q(T))$, $\langle p(T), \delta q(T)\rangle$ can be made arbitrarily large. Thus, by weakening the condition on the variations at the terminal endpoint, we instead obtain the system:
\begin{alignat*}{2}
    \dot{q}(t) &= D_pH(q(t),p(t)) \in \Delta(q(t)),&\,&\text { for all }t\in(0,T), \\
    \dot{p}(t) &+ D_qH(q(t),p(t)) \in \Delta^\circ(q(t)),&\,&\text { for all }t\in(0,T), \\
    q(0) &= q_0, \\
    p(T) &\in \Delta^\circ(q(T)).
\end{alignat*}
This gives an analogue of Type II boundary conditions in the constrained setting. The condition $p(T) \in \Delta^\circ(q(T))$, particularly when the distribution is integrable, can be recognized as a transversality condition; for example, such transversality conditions arise in the Pontryagin maximum principle for a variable endpoint control problem \cite{vp4, Bloch2015}. A related construction is the reduction of Dirac structures, which includes as special cases Lie--Poisson reduction and Euler--Poincar\'{e} reduction, and their associated variational principles \cite{YOSHIMURA2007}. It would be interesting to investigate analogues of Type II variational principles for reduced systems.

\subsection{Stochastic Hamiltonian Systems}\label{sec:stochastic}
We briefly discuss here the Type II variational principle for stochastic Hamiltonian systems. For a thorough discussion of stochastic Hamiltonian systems, see \cite{LaOr2008}. For a recent overview of the symplectic integration of stochastic Hamiltonian systems, see \cite{HoSu2022}. 
Stochastic variational integrators based on a stochastic Hamilton--Pontryagin principle were developed in \cite{BoOw2008}. A Type II variational principle for stochastic Hamiltonian differential equations and variational integrators arising from this variational principle were introduced in \cite{HoTy2018}. This was extended to stochastic Hamiltonian systems with forcing in \cite{KrTy2020}.

Here, we will recall the stochastic Type II variational principle and provide an example of stochastic adjoint processes. Consider a pair of Hamiltonians $H: T^*Q \rightarrow \mathbb{R}$ and $h: T^*Q \rightarrow \mathbb{R}$. Let $(\Omega, \mathcal{F}, P)$ be a probability space with filtration $\{\mathcal{F}_t\}_{t\geq 0}$ and let $W$ be the standard one-dimensional Wiener process. Assuming sufficient regularity on the Hamiltonians \cite{HoTy2018}, a solution curve $(q,p)$ of the stochastic Hamiltonian system
\begin{align*}
    dq &= D_pH dt + D_ph \circ dW(t), \\
    dp &= - D_qH dt - D_qh \circ dW(t),
\end{align*}
where $\circ$ denotes Stratonovich integration, is a critical point of the stochastic Type II action functional $\mathcal{B}: \Omega \times C([0,T]) \rightarrow \mathbb{R}$ under Type II variations. Here, the space of curves is
\begin{align*}
     C([0,T]) := &\{(q,p): \Omega \times [0,T] \rightarrow T^*Q : q,p \text{ are almost} \\
     & \qquad  \text{surely continuous } \mathcal{F}_t-\text{adapted semimartingales}\},
\end{align*}
the action is given by
$$ \mathcal{B}[q,p] := p(T)q(T) - \int_0^T [p \circ dq - H(q,p)dt - h(q,p)\circ dW(t)],$$
and Type II variations are curves $(\delta q,\delta p) \in C([0,T])$ such that almost surely $\delta q(0) = 0$ and $\delta p(T) = 0$. Stochastic discrete Hamiltonian variational integrators are subsequently derived in \cite{HoTy2018} by considering discretizations of the Type II generating function associated with this variational principle.

\begin{example}[Stochastic Adjoint Processes]
    Analogous to the deterministic theory, stochastic adjoint processes arise in optimization and control problems subject to the dynamics of a forward stochastic process (see, e.g., \cite{StochasticAdjoint1, StochasticAdjoint2, StochasticAdjoint3}). 

    Given a forward Stratonovich stochastic differential equation
    \begin{align*}
        dq &= \mu(q)\, dt + \sigma(q) \circ dW(t), \\
        q(0) &= q_0,
    \end{align*}
    define the pair of adjoint Hamiltonians by $H(q,p) = \langle p, \mu(q)\rangle$ and $h(q,p) = \langle p, \sigma(q)\rangle$. The corresponding stochastic Hamiltonian system is given by
    \begin{align*}
        dq &= \mu(q)\, dt + \sigma(q) \circ dW(t), \\
        dp &= - D_q\mu(q)^*p \, dt - D_q\sigma(q)^*p \circ dW(t). 
    \end{align*}
    In the setting of optimization and control problems subject to the dynamics of a forward stochastic process, one considers the Type II free boundary conditions where $q(0)=q_0$ and $p(T) = p_1(q(T))$ for an appropriate choice of mapping $p_1: Q \rightarrow Q^*$. The evolution of $q$ is given by a forward stochastic differential equation, for $s \in [0,T]$,
    $$ q(s) = q_0 + \int_0^s \mu(q(t))\, dt + \int_0^s \sigma(q(t)) \circ dW(t),  $$
    whereas the evolution of $p$ is given by a backward stochastic differential equation (see \cite{StochasticAdjoint4, StochasticAdjoint1}), for $s \in [0,T]$,
    $$ p(s) = p_1(q(T)) + \int_s^T D_q\mu(q(t))^*p(t)\, dt + \int_s^T D_q\sigma(q(t))^*p(t) \circ dW(t).$$
    This forward-backward interpretation of stochastic adjoint processes is naturally compatible with Type II boundary conditions. For examples of applications of stochastic adjoint processes, see, e.g., \cite{StochasticAdjoint1, StochasticAdjoint5, StochasticAdjoint6}.
\end{example}

\section{Conclusion}
We developed several new Type II variational principles for Hamiltonian systems based on a virtual work principle which enforce Type II boundary conditions. Particularly, we show how this construction can be done globally on parallelizable manifolds. By extending the notion of Type II boundary conditions to be free boundary conditions, we also showed how this variational principle can be formulated globally on general manifolds. We further extended this Type II variational principle to the infinite-dimensional setting. Finally, we provided an overview of applications and examples of Hamiltonian variational principles that motivate the utility of studying Hamiltonian systems from the variational perspective.

One natural future research direction is to leverage the intrinsic Type II variational principles introduced in this paper to construct generalizations of Hamiltonian variational integrators to general manifolds, by considering the analogue of Type II generating functions associated with these variational principles. This will have important implications for the intrinsic discretization of optimization algorithms and adjoint systems on manifolds.

\section*{Acknowledgements}
The authors would like to thank the referee for their helpful comments and suggestions. BKT was supported by the Marc Kac Postdoctoral Fellowship at the Center for Nonlinear Studies at Los Alamos National Laboratory. ML was supported in part by NSF under grants DMS-1345013, CCF-2112665, DMS-2307801, and by AFOSR under grant FA9550-23-1-0279. Los Alamos National Laboratory report LA-UR-24-30099.

%\nocite{*}

\bibliographystyle{plainnat}
\bibliography{vp_hamiltonian.bib}

\end{document}